\documentclass[12pt,a4paper]{article}

\usepackage[reqno]{amsmath}
\usepackage{amssymb,amsthm,upref,amscd}
\usepackage[T1]{fontenc}
\usepackage{times}
\usepackage[hidelinks]{hyperref}
\usepackage{color}
\newtheorem{theorem}{Theorem}[section]

\newtheorem{definition}[theorem]{Definition}

\newtheorem{lemma}[theorem]{Lemma}
\newtheorem{proposition}[theorem]{Proposition}
\newtheorem{remark}[theorem]{Remark}

\theoremstyle{definition} \theoremstyle{remark}
\numberwithin{equation}{section}

\def\R{\mathbb{R}}

\newcommand{\de}{\partial}

\newcommand{\D }{\Delta }

\newcommand{\be}{\begin{equation}}
\newcommand{\ee}{\end{equation}}

\begin{document}

\title{Uniqueness of ground state solutions for fractional Schr\"{o}dinger equations with general nonlinearities}
\author{Xinyu Li$^{\mathrm{a}}$\thanks{%
Email: lixinyu@amss.ac.cn.} \ and \ \ Linjie Song$^{\mathrm{b}}$\thanks{% 
Email: songlinjie18@mails.ucas.edu.cn, linjie.song@univ-fcomte.fr.
}\  \\
%EndAName
\\
{\small $^{\mathrm{a}}$ School of Science, Southwest Petroleum University, Chengdu 610500, PR China}\\
{\small $^{\mathrm{b}}$ Universit\'e Marie et Louis Pasteur, CNRS, LmB (UMR 6623), Besan\c{c}on F-25000, France
 }}
\date{}
\maketitle

\begin{abstract}

  In 2013, R.L. Frank and E. Lenzmann [R.L. Frank, E. Lenzmann, Uniqueness of non-linear ground states for fractional Laplacians in $\mathbb{R}$, Acta Math. 210 (2) (2013) 261-318] studied the following problem:
  \begin{align*}
    (-\Delta)^su + u = u^{p-1} \quad \text{in } \R^N,
  \end{align*}
  where $s \in (0,1)$, $N = 1$, $p \in (2,2_s^*)$, and $2_s^*$ is the critical fractional Sobolev exponent. They proved that the ground state is unique (up to translations). Then in 2016, they, together with L. Silvestre [R.L. Frank, E. Lenzmann, L. Silvestre, Uniqueness of radial solutions for the fractional Laplacian, Comm. Pure Appl. Math. 69 (9) (2016) 1671-1726], showed similar uniqueness results for high dimensions ($N \geq 2$), in which they proposed a challenging open problem to extend their results about non-degeneracy and uniqueness of ground states to nonlinearities $f(u)$ beyond the pure-power case. In this paper, we aim to give an affirmative answer to this open issue for two classes of more general nonlinearities, including the case of sums of powers with certain conditions on exponents.
\bigskip



\noindent\textbf{Data availability statement:} Our manuscript has no associated data.
\end{abstract}

%\subjclass{35J65, 58E05}

\newpage
\tableofcontents

\medskip

\section{Introduction and overview of the main results}

Fractional differential operators and the fractional Sobolev spaces attract the 
attention of many researchers, both from a purely mathematical standpoint and 
from a standpoint of various applications, since these operators naturally arise in many different contexts, such as anomalous diffusion and quasi-geostrophic flows, turbulence and water waves, molecular dynamics, and relativistic quantum mechanics of stars (see \cite{BoG,CaV,Co,TZ} and the references therein). These operators also have various applications in probability and finance (see \cite{A,Be,CT}). 
The fractional Schr\"{o}dinger equation discovered by Laskin 
\cite{Laskin1,Laskin2,Laskin3} is a fundamental equation of fractional quantum mechanics. 
Here, what played a central role within these models is often so-called \textit{ground state solutions}, or simply \textit{ground states}. By this, we mean non-trivial, non-negative, and radial functions $u=u(|x|) \ge 0$ that vanish at infinity and satisfy (in the distributional sense) an equation of the form
\begin{align}\label{p1.1}
  (-\Delta)^s u+ F(u) = 0 \quad \text{in } \mathbb{R}^N.
\end{align}
In most examples of interest, the existence of ground states
follows from variational arguments. Based on this variational approach, it is natural in these cases to require that a ground state is also a \textit{minimizer} for some related variational problem. After the existence of ground states is established, some properties, such as the \textit{uniqueness} and the \textit{non-degeneracy} (see Theorem \ref{them2} below for the definition) play an important role in many applications: for example, they are of fundamental importance in the stability and blow-up analysis for the corresponding solitary waves. 
 
However, in striking contrast to the issue of existence, it seems fair to say 
that extremely little was known about the uniqueness of ground states for problems like (\ref{p1.1}), until the outstanding papers \cite{FL13,FS}. For the "classical" limiting case with $s=1$, uniqueness can be addressed by standard ordinary differential equation (ODE) methods. When $s \in (0,1)$, the equation is nonlocal and ODE techniques fail to work. Before \cite{FL13,FS}, the known nonperturbative uniqueness results were for equations with special forms and summarized as follows. (See \cite{FV,KMR} for some perturbative uniqueness results when $s$ is close to 1.)
\begin{itemize}
  \item Ground state solitary waves for the Benjamin-Ono equation in 
  $N=1$ dimension; see \cite{AT}.
  \item Optimizers for fractional Sobolev inequalities in $N \geq 1$ 
  dimensions; see \cite{Li} and \cite{CLO}.
\end{itemize}
In \cite{FL13} and \cite{FS} the  authors studied
\begin{equation}\label{Qp}
(-\Delta)^s u+u-u^{p-1}=0 \quad \text {in } \mathbb{R}^N,
\end{equation}
and proved uniqueness and non-degeneracy of ground states for all $s \in(0,1)$, 
in arbitrary dimension $N \geqslant 1$,  where $2<p<2^*_s$ and the critical fractional Sobolev exponent $2^*_s$ is defined below. In \cite{FS}, the authors proposed a challenging open problem: 
\begin{itemize}
  \item to extend the results about non-degeneracy and uniqueness of ground states obtained in \cite{FL13,FS} to nonlinearities $f(u)$ beyond the pure-power case.
\end{itemize}

In this paper, we aim to give an affirmative answer to this open issue for two classes of more general nonlinearities. More precisely,  we consider the following nonlocal problem:
\begin{equation}\label{Q1}
\left\{
\begin{aligned}
  & (-\Delta)^s u+ \lambda u=f(u) \quad \text{in }\mathbb{R}^N, \\
  & u > 0 \quad  \text{in } \mathbb{R}^N, \quad u(x) \to 0 \quad \text{as } |x| 
  \to \infty,
\end{aligned}
\right.
\end{equation}
where $s\in (0,1)$, $\lambda >0 $ and $ (-\Delta)^s$ denotes the fractional 
Laplacian. The following hypotheses on 
$f:[0,\infty)\to \mathbb{R}$ will be assumed.

\begin{itemize}
  \item[$(f_1)$]  $f \in C^{1,\gamma}((0,\infty), \mathbb{R}) \cap C^{1}([0,\infty), \mathbb{R})$ for some $\gamma>\max \{0,1-2 s\}$, $f(0)=f'(0)=0$, $f(t)\leqslant C(t+t^p)$ for some $p \in(1,2_s^*-1)$ and $C>0$.
  \item[$(f_2)$]  $f' (t)\geqslant 0,$   $f' (t)t>f(t)$ and $f'$ is nondecreasing on  $(0,\infty)$.
\end{itemize}
A classical example satisfying $(f_1)-(f_2)$ is the sum of pure-power nonlinearities:
\begin{itemize}
	\item[$(F)$]  
	$f(t)=\sum_{j=1}^{m}a_jt^{p_j},
	\qquad a_j>0,\qquad1<p_1<\cdots<p_m<p_s,$ where $p_s=2_s^*-1$.
\end{itemize}

Our main results on this class of nonlinearities are as follows.

\begin{theorem}\label{thm}
Let $N\ge1$ and $0<s<1$. Fix a compact interval $J\Subset(1,p_s)$.
Assume that $(F)$ holds and $p_j\in J.$
Suppose that
\[
p_m-p_1\le(p_1-1)C_J^{-(p_m-1)/2}
+(p_m-1)C_J^{-(p_1-1)/2},
\]
where $C_J>1$ is a constant depending only on $N,s,J$, given in Proposition \ref{prop:peak}.
Then the positive ground state solution $u>0$ of (\ref{Q1}) is unique (up to translations) and non-degenerate.
\end{theorem}

\begin{remark}
	When $N\geq 2$, $0<s<1$, and $b<N/(N-2s)$, where $b = \max_{p \in J} p$, we provide an explicit choice of $C_J$ in Proposition~\ref{prop:peak-explicit}.
\end{remark}

\begin{theorem}\label{thmN1}
Let $N=1$. Assume that $(F)$ holds. Then the positive ground state solution $u>0$ of (\ref{Q1}) is unique (up to translations) and  non-degenerate if one of the following assertions holds:
\begin{itemize}
	\item[(i)] 
	\[
	\frac12\le s<1,\qquad1<p_1<\cdots<p_m\le3;
	\]
	\item[(ii)]
	\[
	s_0<s<\frac12, \qquad 
	1<p_1<\cdots<p_m\le\frac4{C_s}-1,
	\]
	where
	\[
	C_s=\frac{\Gamma(\frac32-s)\Gamma(s)}
	{\Gamma(1-s)\Gamma(s+\frac12)},
	\]
	and $s_0\in(0,\tfrac12)$ is the unique solution of $C_{s_0}=2$;
	\item[(iii)]
	\[
	\tfrac12<s<1,\qquad  p_m-p_1\le(p_1-1)\sqrt{K_{s,p_m}}
	+(p_m-1)\sqrt{K_{s,p_1}},
	\]
	where $K_{s,r}, r \ge 1$ is given in \eqref{def:K} later.
	
\end{itemize}
\end{theorem}

The following result addresses another class of nonlinearities.

\begin{theorem}\label{thm:log}
Let $N\ge1$, $0<s<1$, and
\[
f(t)=ct^p\log(1+t^\delta),\qquad c>0,\quad p\ge1,\quad\delta>0.
\]
Assume $p+\delta<p_s$. Then the positive ground state solution $u>0$ of (\ref{Q1}) is unique (up to translations) and non-degenerate.
\end{theorem}

We remark that our results also hold in the case $s=1$, and cover some new results which were not obtained in the literature. In the local case when $s=1$, the uniqueness of the positive radial solution for problem (\ref{Qp}) with all $p \in (2,2^*)$ was established in a celebrated paper by Kwong 
\cite{Kwong} (see also \cite{32}). Things become difficult for a more complicated case with 
double-power nonlinearities, i.e., $f(u) = u^{p-1}+u^{q-1}$ with $2 < p < q < 2^*$. 
Indeed, in $N \geqslant 3$ dimensions, Chen and Lin \cite{CCC1991} established the uniqueness of 
the positive radial solution covering examples $f(u) = u^{p-1}+u^{q-1}$ with 
$2 < p < q \leqslant \frac{2N-2}{N-2}$. Note that there is a gap between 
the exponent $\frac{2N-2}{N-2}$ in \cite{CCC1991} and the critical exponent 
$2^*=\frac{2N}{N-2}$. Later, D\'avila, del Pino and Guerra \cite{DDG} showed that this gap cannot be completely bridged when $N = 3$. Recently, Song and Zou \cite{SZ} proved that the gap can be closed for $N \ge 5$. The result presented in Theorem \ref{thm} involves sum of three or more powers and is new even when $s=1$. 

For more results on uniqueness problems involving a nonlocal operator, readers are referred to \cite{APZ,CGWY,ad,A. Dieb2023,weth} and references therein.
\bigskip

{\bf Sketch of Main Ideas. }
We work with variational ground states, whose Morse index is one.
To complete our proof, we need to overcome two main challenges caused by the nonlinear terms. The first one is to obtain the non-degeneracy of the positive 
radial solution $u$ having Morse index one, i.e., to prove that the kernel of the self-adjoint operator 
$L_{+} :=(-\Delta)^s+\lambda- f'(u)$ acting on $L^2(\R^N)$ is exactly 
$\operatorname{span}\{\de_{x_i}u: i=1, \cdots, N\}$. Inspired by \cite{FL13, FS}, we introduce a function 
\[
H_{f,\lambda}(t) = \frac{f(t)-\lambda t}{tf'(t)-f(t)}
\]
and aim to find conditions to prove that 
\[
H_{f,\lambda}'(t) > 0, \qquad t \in (0,|u|_{\infty}).
\]
Together with the oscillation result in \cite{FS} (see Proposition \ref{prop2.7}), we will show that the kernel of $L_+$ is trivial in the radial space. Then, using spherical harmonics and decomposing the action of $L_{+}$ yields the non-degeneracy in the whole space.  

After establishing the non-degeneracy, we use a continuation argument to prove 
the uniqueness. In contrast to the discussions in \cite{FL13, FS}, the more general nonlinearities bring new difficulties. Here, we use a homotopic argument by continuation in $\lambda$. More precisely, starting from some solution $u_0$ of
$$
(-\Delta)^{s} u_0+\lambda_0 u_0=f(u_0)  \quad \text {in  } \mathbb{R}^N,
$$
we seek to construct a branch $u_\lambda$ with $\lambda$ stemming from 
$\lambda_0$ and solving the problem
$$
(-\Delta)^{s} u+\lambda u=f(u) \quad \text {in  }  \mathbb{R}^N .
$$
The local existence and uniqueness (in a suitable function space) 
for $u_\lambda$ with $u_{\lambda_0}=u_0$ follow from an implicit function 
argument, provided the linearization around $u_0$ is invertible, which is the 
first key step in the argument. The second key point is now to derive suitable a priori bounds that guarantee that the branch $u_\lambda$ can be extended to 
$\lambda \nearrow +\infty$. Finally, we aim to prove that the positive solution whose Morse index is one is 
unique when $\lambda$ is large enough. To do it, we use scaling arguments. To prove Theorem \ref{thm:log}, we instead fix both $s$ and $\lambda$ and deform the nonlinearity by a parameter $\theta\in[0,1]$.

\bigskip

{\bf Plan of the Paper.}
We organize this paper as follows. In Section \ref{Preliminary} we give some preliminary results. 
 Section \ref{nondegeneracy} is devoted to the proof of Theorem \ref{them2} (non-degeneracy result). The proofs of Theorems \ref{thm}, \ref{thmN1} will be given in
Section \ref{uniqueness} and Theorem \ref{thm:log} will be proved separately in Section \ref{sec:log-cont}.  

\bigskip

{\bf Notation and Conventions.} 
Throughout this paper, we employ the common abuse of notation by writing both 
$f=f(|x|)$ and $f=f(r)$ for any radial function $f$ on $\mathbb{R}^N$. We use 
standard notation for $L^p$ and Sobolev spaces and 
$L_{\mathrm{rad}}^p(\mathbb{R}^N)$ denotes the space of radial and 
$p$-integrable functions on $\mathbb{R}^N$. For $k \in \mathbb{N}$ and 
$0<\gamma \leqslant 1$ the Hölder space $C^{k, \gamma}(\mathbb{R}^N)$ 
is equipped with the norm $\|u\|_{C^{k, \gamma}}=\sum_{|\alpha| \leqslant k}\left\|\partial_x^\alpha u\right\|_{L^{\infty}}+\sum_{|\alpha|=k}[\partial_x^\alpha u]_{C^\gamma}$, where $[u]_{C^\gamma}=\sup _{x \neq y} \frac{|u(x)-u(y)|}{|x-y|^\gamma}$. We often write $L^p$ instead of $L^p(\mathbb{R}^N)$ etc.  We denote  $ \langle \cdot, \cdot \rangle_{L^2}$  and  $|\cdot|_2$  the inner product  and  norm  on $L^2(\mathbb{R}^N)$, respectively. $B_R(x_0)=\{x:|x-x_0|<R\}$ and $B_R=B_R(0)$.
We denote $|B_R|$ for volume and $|\mathbb S^{N-1}|$  surface measure.
 The
Fourier transform is defined by
$\widehat v(\xi)=\int_{\mathbb R^N}e^{-ix\cdot\xi}v(x)\,dx$;
accordingly $\widehat{(-\Delta)^sv}(\xi)=|\xi|^{2s}\widehat v(\xi)$.
The letter $\Gamma$ denotes the gamma function, and
$\mathrm B(a,b)=\Gamma(a)\Gamma(b)/\Gamma(a+b)$  the beta function.
The fractional critical power is
$p_s=(N+2s)/(N-2s)$ for $N>2s$, and $p_s=+\infty$ otherwise;
the corresponding critical Sobolev exponent is
$2_s^*=2N/(N-2s)$ for $N>2s$, and $2_s^*=+\infty$ otherwise.

We employ the following convention for constants in this paper: Unless otherwise explicitly stated, we write
$$
X \lesssim_{a, b, c, \ldots} Y
$$
to denote that $X \leqslant C Y$ with some constant $C>0$ that only depends on 
the quantities $a, b, c, \ldots$ and the space dimension $N \geqslant 1$. 
 Moreover, we require that $C>0$ can be chosen uniform if $a, b, c, \ldots$ 
 range in some fixed compact set. We write $X \lesssim Y$ for simplicity without 
 causing ambiguity. We denote by $X \sim Y$ if $X \lesssim Y$ and $Y \lesssim X$.

\bigskip

\section{Preliminaries}\label{Preliminary}

Let $s\in (0,1)$ and   $ (-\Delta)^s$ denote the fractional Laplacian given by
\[
(-\Delta)^su(x)=c_{N,s} \ P.V. \int_{\mathbb{R}^N}  \frac {u(x)-u(y)}{|x-y|^{N+2s}}dy, 
\quad  c_{N,s} := 4^{s}\pi^{-\frac N2} \frac{\Gamma(\frac N2+s)}{\Gamma(2-s)}s(1-s).
\] We define $H^{s}(\mathbb{R}^N)$ the fractional Sobolev space,
\begin{align}
H^s(\mathbb{R}^N)=\left\lbrace u \in
L^2(\mathbb{R}^N):\frac{|u(x)-u(y)|}{|x-y|^{\frac{N}{2}+s}} \in
L^2(\mathbb{R}^N\times\mathbb{R}^N)\right\rbrace ,
\end{align}
endowed with the norm
\begin{align}\label{norm:def}
\|u\|:= \sqrt{|u|_2^2+\mathcal{E}(u,u)},
\end{align}
where, for $u,v\in H^s(\mathbb{R}^N)$,
\begin{align*}
\mathcal{E}(u,v):=\dfrac{c_{N,s}}{2}\iint_{\mathbb{R}^N\times \mathbb{R}^N}
\dfrac{(u(x)-u(y))(v(x)-v(y))}{|x-y|^{N+2s}}dx\,dy.
\end{align*}

We briefly recall the extension principle that expresses the fractional Laplacian as a Dirichlet-Neumann map for a local elliptic problem in one additional variable; see \cite[Section~4.1]{FS} and \cite{CSext}.
For a measurable $u : \mathbb R^N \to \mathbb R$, we define its {\em $s$-Poisson extension} to the upper half-space
$ \mathbb R^{N+1}_+
 =\{(x,y):x\in\mathbb R^N,\ y>0\}$
by 
\[
 U(x,y)=\int_{\mathbb R^N}
 P^{\mathrm{ext}}_{N,s,y}(x-z)u(z)\,dz,
\]
where the generalized Poisson kernel is
\[
 P^{\mathrm{ext}}_{N,s,y}(x)
 =\frac{\Gamma(\frac N2+s)}{\pi^{N/2}\Gamma(s)}
   \frac{y^{2s}}{(|x|^2+y^2)^{(N+2s)/2}}.
\]
This kernel is positive and satisfies
$\int_{\mathbb R^N}P^{\mathrm{ext}}_{N,s,y}(x)\,dx=1.$
Under suitable assumptions on $u$, the extension $U$ solves the following degenerate elliptic boundary value problem
\[
 \begin{cases}
  \operatorname{div}(y^{1-2s}\nabla U)=0
   &\text{in }\mathbb R^{N+1}_+,\\
  U=u
   &\text{on }\partial \mathbb R^{N+1}_+,
 \end{cases}
\]
where the boundary condition is understood in the trace sense. Moreover, if $u$ is sufficiently regular, there holds
\[
 (-\Delta)^s u
 =-d_s^{\mathrm{ext}}\lim_{y\downarrow0}
 y^{1-2s}\partial_yU(\cdot,y),
 \quad \text{where}\quad
 d_s^{\mathrm{ext}}=\frac{2^{2s-1}\Gamma(s)}{\Gamma(1-s)}.
\]
This identity gives the Dirichlet--Neumann representation with the
weight $y^{1-2s}$. When $s=\tfrac12$, the extension equation reduces to
the ordinary Laplace equation.

For $N=1$, we write $P_{s,y}=P^{\mathrm{ext}}_{1,s,y}$.
We set
\[
 V(x,y)=-y^{1-2s}\partial_yU(x,y).
\]
It satisfies
\[
 \operatorname{div}(y^{2s-1}\nabla V)=0,
 \qquad
 V(\cdot,0)=\frac{(-\Delta)^su}{d_s^{\mathrm{ext}}}.
\]
Therefore $V$ satisfies the extension equation of order $1-s$; see
\cite[Section~2.3]{CSext} and \cite[Proposition~3.6]{CaSi}.

\medskip

Let $P_s(\tau,x)$ denote the fractional heat kernel associated with
$e^{-\tau(-\Delta)^s}$. It is characterized by
\[
 \widehat{P_s(\tau,\cdot)}(\xi)=e^{-\tau|\xi|^{2s}},
 \qquad \tau>0.
\]
For initial data $v$, we have that
\[
 W(\tau,\cdot)=P_s(\tau,\cdot)*v
\]
solves
\[
 \partial_\tau W+(-\Delta)^sW=0,
 \quad \text{where}\quad W(0,\cdot)=v.
\]
The kernel is positive, has integral one, and satisfies 
\[
 P_s(\tau,x)
 =\tau^{-N/(2s)}P_s\bigl(1,\tau^{-1/(2s)}x\bigr).
\]
Here $\tau$ is the time variable, while $y$ is the extension variable.
The positivity and scaling of $P_s$ will be used in the blow-up argument;
see also \cite[Appendix~C.1]{FS}.

\medskip

For $\lambda>0$, we denote by $G_{s,\lambda}$ the kernel of the resolvent
$\bigl(( -\Delta)^s+\lambda\bigr)^{-1}$. Then
\[
 \widehat G_{s,\lambda}(\xi)=\frac1{|\xi|^{2s}+\lambda},
 \qquad
 G_{s,\lambda}(x)=\int_0^\infty
 e^{-\lambda\tau}P_s(\tau,x)\,d\tau.
\]
This kernel is positive and integrable, with
\[
 \int_{\mathbb R^N}G_{s,\lambda}(x)\,dx=\frac1\lambda.
\]
Consequently, solutions of
 $(-\Delta)^su+\lambda u=f(u)$
can be written as
\[
 u=G_{s,\lambda}*f(u).
\]
In dimension one, we write $G_s=G_{s,1}$.
The value of $G_s$ at the origin is finite only if $s>\tfrac12$, in which case
\[
 c_s:=G_s(0)=\frac1{2s\sin(\pi/(2s))},\quad 
 h_s(x):=\frac{G_s(x)}{c_s}.
\]
and
\begin{align} \label{def:K}
	K_{s,r}:=\int_{\mathbb R}G_s(x)h_s(x)^r\,dx\quad(r\ge1).
\end{align}

\smallskip

\subsection{Properties of the positive  solutions}
 
In this paper, we focus on real-valued solutions. Throughout the paper, we extend $f$ by zero to $(-\infty,0)$. Since $f(0)=f'(0)=0$, this extension is of class $C^1$. A natural approach to construct 
nontrivial solutions $u$ for equation \eqref{Q1} is given by minimizing $J$ on 
the Nehari manifold  $\mathcal{N}$, where $J: H^s(\mathbb{R}^N) \rightarrow \mathbb{R}$ is the energy functional given by
$$
J(v):=\frac{1}{2} \mathcal{E}(v, v)+\frac{\lambda}{2}|v|_{2}^2- 
\int _{\mathbb{R}^N}F(v)dx, \quad F(t) := \int_0^tf(s)ds,
$$
and the set $\mathcal{N}$ is the Nehari manifold

$$
\mathcal{N}:=\left\{ v \in H^s (\mathbb{R}^N ) \backslash\{0\}: \mathcal{E}(v, v)+\lambda|v|_{2}^2-\int _{\mathbb{R}^N}f(v)vdx =0 \right\}.
$$

We summarize the following  fundamental properties 
of positive solutions for equation \eqref{Q1}.

\begin{proposition} \label{rem3.1}
Let $N \geqslant 1, s \in(0,1)$, $\lambda >0$ and  assume that $(f_1)$ 
holds true. Then the following  holds.\\
{\bf  Symmetry, regularity, and decay }: If $u \in H^s(\mathbb{R}^N)$ with 
$u > 0$ solves  \eqref{Q1}, then there exists some $x_0 \in \mathbb{R}^N$ such 
that $u$ is radial about $x_0$, positive, and strictly decreasing in 
$\left|x-x_0\right|$. Moreover, the function $u$ belongs to 
$H^{2 s+1}(\R^N) \cap C^{2s+\alpha}(\mathbb{R}^N)$ for some $\alpha > 0$
and it satisfies
$$
\frac{C_1}{1+|x|^{N+2 s}} \leqslant u(x) \leqslant \frac{C_2}{1+|x|^{N+2 s}} \quad 
\text {  for } x \in \mathbb{R}^N,
$$
with some constants $C_2 \geqslant C_1>0$ depending on $s, N, \lambda$, and $u$.
\end{proposition}

\begin{remark}
  See  \cite[Theorems 1.1-1.3 ]{Byeon} for the existence result and the results of  symmetry, regularity, and decay in Proposition \ref{rem3.1}.
  As for the strictly decreasing results, we can apply the moving plane method for the fractional 
  Laplacian in \cite{CWX2017} to obtain that $u(x)$ is radial about $x_0$, positive, and (but not necessarily strictly) decreasing in $\left|x-x_0\right|$. To show that $u(x)$ is indeed strictly decreasing in $\left|x-x_0\right|$, we use an adaptation of the argument
in \cite[ Appendix D]{FS}.
\end{remark}

\subsection{Definitions of the non-degeneracy and Morse index}
The solution obtained by  minimizing $J$ on 
the Nehari manifold  has more properties, for example, its Morse index is no more than $1$. With additional hypotheses stated below, the Morse index is exactly $1$ and the solution is also non-degenerate. For further explanations, 
we introduce the following notions.

\begin{definition}[Morse index] \label{defmi}
  Let $u \in H^s(\mathbb{R}^N)$ be a solution of 
  (\ref{Q1}), and we define the linearized operator
$$
L_{+}=(-\Delta)^s+\lambda- f'(u).
$$
Then, $L_{+}$ is a self-adjoint operator on $L^2(\mathbb{R}^N)$ with quadratic-form domain $H^s(\mathbb{R}^N)$ and operator domain $H^{2 s}(\mathbb{R}^N)$. The Morse index of $u$ is defined as
$$
\mu\left(u \right):=\sharp\left\{e<0: e \text { is an eigenvalue of } L_{+}\right\},
$$
where multiplicities of eigenvalues are taken into account. Noticing that $L_{\text {rad }}^2(\mathbb{R}^N)$ is an invariant subspace of $L_{+}$, we also introduce the Morse index of $u$ in the sector of radial functions by defining
$$
\mu_{\text {rad }}\left(u\right):=\sharp\left\{e<0: e \text { is an eigenvalue of } L_{+} \text { restricted to } L_{\text {rad }}^2(\mathbb{R}^N)\right\} \text {. }
$$
\end{definition}

\begin{remark} \label{rmknon}
If $u$ is the solution of (\ref{Q1}) such 
that $J(u) = \inf_{\mathcal{N}}J(v)$, and assume  $(f_2)$ holds, then $\mu\left(u\right)=\mu_{\text {rad }}\left(u\right)=1$.
In fact, on the one hand, by $(f_2)$,
$$
\left\langle L_{+} u, u\right\rangle_{L^2}= \int_{\mathbb{R}^N}  f(u)u-f'(u) u^2 dx<0.
$$
Thus, by the min-max principle, the operator $L_{+}$ has at least one negative 
eigenvalue, i.e. $\mu\left(u\right) \geqslant 1$. Since $u$ is radial, then 
$\mu_{\text {rad }}\left(u\right) \geqslant 1$ by applying min-max principle in 
$L_{r a d}^2(\mathbb{R}^N)$. On the other hand, since $u$ is a minimum point 
in the Nehari manifold whose codimension is $1$, it is easy to verify that 
$\mu_{\text {rad }}\left(u\right) \leqslant$ $\mu\left(u\right) \leqslant 1$.
\end{remark}

Now we state our result on the non-degeneracy, whose definition will also be 
given.

\begin{theorem}[Non-degeneracy] \label{them2} 
  Let $u \in H^s(\mathbb{R}^N)$ be a positive solution of 
  (\ref{Q1}) having Morse index $1$. Under the conditions of Proposition \ref{rad}, then $u$ is non-degenerate, in the sense 
  that 
$$
\ker L_{+}=\operatorname{span}\left\{\partial_{x_1} u, \ldots, \partial_{x_N} u
\right\},
$$
where $L_+$ is given in Definition \ref{defmi}.
\end{theorem}

\subsection{An oscillation estimate of the second radial eigenfunction }

In this subsection,  we collect some preliminary information, which is useful to 
prove Theorem \ref{them2}. We start with the following result from \cite{FS}.

Let $N \geqslant 1$ and $s \in(0,1)$ be given. We consider the linear equation
\begin{equation}\label{2.1V}
(-\Delta)^s u+V u=0 \quad \text{in } \mathbb{R}^N.
\end{equation}
We require that the potential $V: \mathbb{R}^N \rightarrow \mathbb{R}$ satisfies the following conditions:
\begin{itemize}
  \item[(V1)] $V=V(|x|)$ is radial and non-decreasing in $|x|$.
  \item[(V2)] $V \in C^{0, \gamma}(\mathbb{R}^N)$ for some $\gamma>\max \{0,1-2 s\}$. 
\end{itemize}

\begin{proposition}[Frank-Lenzmann-Silvestre, \cite{FS}, Corollary 1] \label{prop2.5}  
Suppose $N \geqslant 1, s \in(0,1)$, and let $V$ satisfy (V1) and (V2). 
Consider $H=(-\Delta)^s+V$ acting on $L_{\mathrm{rad}}^2(\mathbb{R}^N)$. Then all eigenvalues of $H$ are simple.

In particular, if $E_1^{\mathrm{rad}} \leqslant E_2^{\mathrm{rad}} \leqslant E_3^{\mathrm{rad}} \leqslant \cdots<\inf \sigma_{\mathrm{ess}}(H)$ denote (with counting multiplicity) the discrete eigenvalues of $H$ acting on $L_{\mathrm{rad}}^2(\mathbb{R}^N)$, then we have strict inequalities
$$
E_1^{\mathrm{rad}}<E_2^{\mathrm{rad}}<E_3^{\mathrm{rad}}<\cdots<\inf \sigma_{\mathrm{ess}}(H).
$$
\end{proposition}

As the last main result for the linear equation (\ref{2.1V}), we need the following sharp oscillation estimate for the second eigenfunction of $H$, which provides us with a bound in agreement with Sturm-Liouville theory for ODEs.

\begin{proposition}[Frank-Lenzmann-Silvestre, \cite{FS}, Theorem 2] \label{prop2.7}
Let $H=(-\Delta)^s+V$ be as the one given by Proposition \ref{prop2.5} above. 
Suppose that $H$ has at least two radial eigenvalues $E_1^{\mathrm{rad}}<
E_2^{\mathrm{rad}}<\inf \sigma_{\mathrm{ess}}(H)$. Let $\psi \in 
L^2(\mathbb{R}^N)$ denote the radial eigenfunction of $H$ for the second radial 
eigenvalue $E_2^{\mathrm{rad}}$. Then $\psi=\psi(|x|)$ changes its sign exactly 
once for $|x|=r \in(0,+\infty)$. By this, we mean that there is some $\rho>0$ 
such that (after multiplying $\psi$ with $-1$ if necessary) we have
$$
\psi(r) \leqslant 0 \text { for } r \in\left[0, \rho\right) \text { and } 
\psi(r) \geqslant 0 \text { for } r \in\left(\rho,+\infty\right),
$$
and $\psi \not \equiv 0$ on both intervals $\left[0, \rho\right)$ and 
$\left[\rho,+\infty\right)$. Note also that $\psi(0)<0$.
\end{proposition}

\begin{remark}
	For the applications in the following sections, the potential $V$ is $\lambda - f'(u)$, where $u$ is a positive solution of \eqref{Q1} having Morse index $1$. Condition $(V1)$ follows from the decrease of $u$ and the monotonicity of $f'$ already included in $(f_2)$. Condition $(V2)$ additionally requires the stated regularity assumptions. 
\end{remark}

\section{Non-degeneracy of the ground state solution}\label{nondegeneracy}

This section is devoted to the proof of Theorem \ref{them2} that establishes the non-degeneracy of solutions $u > 0$ with Morse index $1$ for equation (\ref{Q1}). By Proposition 
\ref{rem3.1}, we can assume that $u=u(|x|)>0$ is radial without loss of 
generality.

\subsection{Non-degeneracy in the non-radial sector}

In the first step, we show that the kernel of the linearized operator $L_{+}$ 
in the non-radial sector is exactly $\operatorname{span} \{\de_{x_i}u: i=1, \cdots, N\}$.

For $N\ge2$, since $u=u(|x|)$ is a radial function, the operator $L_{+}=(-\Delta)^s +\lambda 
- f'(u) $ commutes with rotations on $\mathbb{R}^N$.   Using the decomposition in terms of spherical harmonics

\begin{equation}\label{s7.5}
L^2(\mathbb{R}^N)=\bigoplus_{\ell \geqslant 0} \mathcal{H}_{\ell},
\end{equation}
we find that $L_+$ acts invariantly on each subspace
\begin{equation}\label{s7.6}
\mathcal{H}_{\ell}=L^2\left(\mathbb{R}^+, r^{N-1} d r\right) \otimes \mathcal{Y}_{\ell} .
\end{equation}
Here $\mathcal{Y}_{\ell}=\operatorname{span}\left\{Y_{\ell, m}\right\}_{m \in M_{\ell}}$ denotes space of the spherical harmonics of degree $\ell$ in space dimension $N$. Note that the index set $M_{\ell}$ depends on $\ell$ and $N$. Recall also that $-\Delta_{\mathbb{S}^{N-1}} Y_{\ell, m}=\ell(\ell+N-2) Y_{\ell, m}$.
For each $\ell \geqslant 0$, the action of $L_{+}$ on the radial factor in $\mathcal{H}_{\ell}$ is given by
\begin{align}\label{u7.7}
\left(L_{+, \ell} g\right)(r)=\left(\left(-\Delta_{\ell}\right)^s g\right)(r)+\lambda 
g(r)- f'(u)g(r),
\end{align}
for $g\in C_0^{\infty}(\mathbb{R}^+)\subset L^2(\mathbb{R}^+, r^{N-1}dr)$. Here $\left(-\Delta_{\ell}\right)^s$ is given by spectral calculus and the known formula
\begin{align}\label{u7.8}
-\Delta_{\ell}=-\frac{\partial^2}{\partial r^2}-\frac{N-1}{r} \frac{\partial}{\partial r}+\frac{\ell(\ell+N-2)}{r^2} .
\end{align}
Technically, we consider (\ref{u7.8}) as a self-adjoint operator in 
$L^2(\mathbb{R}^+,r^{N-1}dr)$ defined as the Friedrichs extension of the corresponding 
differential expression acting on $C_0^{\infty}(\mathbb{R}^+)$.

We now derive the following result inspired by \cite{FS}.
\begin{lemma}\label{nonrad}  Let $N \geq 1,0<s<1$, and $\lambda>0$. Suppose that $u$ is a positive, radially symmetric, and strictly radially decreasing solution of \eqref{Q1}, with sufficient regularity and decay. Then

$$
\operatorname{ker} L_{+} \cap L_{\mathrm{rad}}^2\left(\mathbb{R}^N\right)^{\perp}=\operatorname{span}\left\{\partial_{x_1} u, \ldots, \partial_{x_N} u\right\} .
$$

\end{lemma}

\begin{proof}

Suppose $N\ge2$ and $u$ has the regularity, decay and strict radial monotonicity stated in Proposition~\ref{rem3.1}. Claim
\[
\ker L_{+,1}=\operatorname{span}\{\partial_ru\},\qquad
\ker L_{+,\ell}=\{0\}\quad(\ell\ge2).
\]

By differentiating equation \eqref{Q1}, we see that $L_+ \partial_{x_i}u=0$ for $i=1,\cdots ,N$, 
where $\partial_{x_i}u=\frac{\partial u}{\partial x_i}=u'(r)\frac{x_i}r \in \mathcal{H}_{1}$.
Hence, $L_{+,1}u'(r) = 0$. By standard arguments (because of $u'<0$ and the Perron-Frobenius 
property of $L_{+,1}$), we deduce that $0$ is the lowest and hence non-degenerate eigenvalue of 
$L_{+,1}$. Thus we conclude that
\[
\operatorname{ker} L_{+,1}=\operatorname{span}\{\partial_r u\}.
\]

Moreover, by \cite[Lemma C.5]{FS}, we have the strict inequality $L_{+,\ell}>L_{+,1} $ in the sense of quadratic forms for any $\ell\geqslant 2$.
Since $0$ is the lowest eigenvalue of $L_{+,1}, $ this implies that $0$ cannot be an eigenvalue of $L_{+,\ell} $ for  any $\ell\geqslant 2$.

For $N=1$, the angular decomposition is replaced by the even-odd decomposition. Since $u$ is even, $L_+$ acts invariantly on each subspace. Differentiating the equation yields $L_{+} u^{\prime}=0$. Moreover, $u^{\prime}(x)<0$ for $x>0$. By the Perron-Frobenius property (see \cite[Theorem XIII.44]{RS4}) in the odd subspace (see \cite[Lemma C.3]{FL13}), $0$ is the simple lowest eigenvalue of the restriction of $L_{+}$ to odd functions. Consequently,

$$
\operatorname{ker}\left(\left.L_{+}\right|_{L_{\text {odd }}^2}\right)=\operatorname{span}\left\{u^{\prime}\right\} .
$$
This completes the proof.

\end{proof}

\subsection{ Non-degeneracy in the radial sector}

In this subsection, we prove that the kernel of $L_+$ is trivial in the radial sector.

\begin{proposition}\label{rad}  Let \(N\geq 1\), \(0<s<1\), and \(\lambda>0\). Suppose that \(u\in H^s(\mathbb R^N)\) is a positive, radial, strictly decreasing solution of  (\ref{Q1})
with Morse index $1$. 
Assume that
\[
H_{f, \lambda}(t):=\frac{f(t)-\lambda t}{Q(t)}
\]
is strictly increasing on \((0,u(0))\), where $Q(t)=tf'(t)-f(t)$, and $f$ satisfies $(f_1)$ and $(f_2)$. Then
\[
\ker L_+\cap L_{\mathrm{rad}}^2(\mathbb R^N)=\{0\}.
\]

\end{proposition}

\begin{proof} Suppose by contradiction that there exists $w \in\left(\operatorname{ker} L_{+}\right) \cap L_{\text {rad }}^2$ and $w \neq 0$. Since the Morse index of $u$ is $1$ , then $0$ is the second eigenvalue of $L_{+}$ with eigenfunction $w$. From Proposition~\ref{prop2.7}, we conclude that there is some $\rho>0$ such that (after multiplying $w$ with $-1$ if necessary)

\[
w(r)\geq0\quad(0<r<\rho),
\qquad
w(r)\leq0\quad(r>\rho),
\]
and $w \not \equiv 0$ on both intervals $[0, \rho)$ and $[\rho,+\infty)$.

Note that a calculation shows that
$L_+u=-Q(u)$ and $L_+(x\cdot\nabla u)=2s(f(u)-\lambda u).$
By self-adjointness of \(L_+\) and \(L_+w=0\), it follows that
\[
\int_{\mathbb R^N}Q(u)w\,dx=0,
\qquad
\int_{\mathbb R^N}(f(u)-\lambda u)w\,dx=0.
\]

Set \(T=u(\rho)\). Combining these two identities yields
\[
\begin{aligned}
0
&=\int_{\mathbb R^N}
Q(u)\bigl[H_{f, \lambda}(u)-H_{f, \lambda}(T)\bigr]w\,dx\\
&=|\mathbb S^{N-1}|
\int_0^\infty
Q(u(r))\bigl[H_{f, \lambda}(u(r))-H_{f, \lambda}(T)\bigr]
w(r)r^{N-1}\,dr.
\end{aligned}
\]

However, \(u\) is strictly decreasing and \(H_{f, \lambda}\) is strictly increasing. Hence, we have
\[
H_{f, \lambda}(u(r))-H_{f, \lambda}(T)
\begin{cases}
>0,&0<r<\rho,\\
<0,&r>\rho.
\end{cases}
\]
Since \(Q(t)>0\) for every \(t>0\), the integrand in the last identity is nonnegative everywhere and strictly positive on a set of positive measure. Consequently,
\[
\int_0^\infty
Q(u(r))\bigl[H_{f, \lambda}(u(r))-H_{f, \lambda}(T)\bigr]
w(r)r^{N-1}\,dr>0,
\]
which is a contradiction. Therefore, we have
\[
\ker L_+\cap L_{\mathrm{rad}}^2(\mathbb R^N)=\{0\}.
\]
The proof is complete.
\end{proof}

\begin{remark} The pure-power nonlinearity $f(t)=t^p$, $1<p<2_s^*-1$ is a special case of this scalar monotonicity condition. Indeed, for $p>1$,
$Q(t)=tf'(t)-f(t)=(p-1)t^p,$
and hence
$H_{f,\lambda}(t)
=\frac{1}{p-1}(1-\lambda t^{1-p}).$
Therefore,
$H_{f,\lambda}'(t)=\lambda t^{-p}>0$
for every $t, \lambda>0.$
\end{remark}

Now we can complete the proof of Theorem \ref{them2}. 
\begin{proof}[Proof of Theorem \ref{them2}]
	Proposition~\ref{rad} gives 
	$(\operatorname{ker} L_+)\cap  L^2_{rad}(\mathbb{R}^N)=\{0\}$. Together 
	with Lemma \ref{nonrad}, we conclude that 
	$$
	\ker L_{+}=\operatorname{span}\left\{\partial_{x_1} u, \ldots, \partial_{x_N} u
	\right\}.
	$$
\end{proof}

In the  following subsections, we search for assumptions on $f$ such that $H_{f,\lambda}$ is strictly increasing for $0 < t < |u|_\infty$.

\subsection{The case of sums of powers: $N \ge 1$}

\begin{lemma}\label{sufficient condition} Let
$
f(t)=\sum_{j=1}^m a_j t^{p_j},  a_j>0,  1<p_1<\cdots<p_m,
$
and  $\lambda>0$. For a fixed $t>0$, suppose that there exist positive numbers $k_1, \ldots, k_m$ and a real number $\zeta$ such that
\begin{align} \label{eqkzeta}
	\sum_{j=1}^m a_j k_j t^{p_j-1}<\lambda; \qquad \left[1-\zeta\left(p_j-1\right)\right]^2 \leq k_j, ~~ j=1, \ldots, m.
\end{align}
Then $
H_{f, \lambda}^{\prime}(t)>0 .
$
\end{lemma}

\begin{proof} Since
$$
H_{f, \lambda}(t)=\frac{f(t)-\lambda t}{t f^{\prime}(t)-f(t)}=\frac{\sum_{j=1}^m a_j t^{p_j-1}-\lambda}{\sum_{j=1}^m a_j\left(p_j-1\right) t^{p_j-1}},
$$
we have
\begin{equation}\label{eq:moment-derivative}
t H_{f, \lambda}^{\prime}(t)=\frac{\left(\sum_{j=1}^m a_j\left(p_j-1\right) t^{p_j-1}\right)^2-\left(\sum_{j=1}^m a_j t^{p_j-1}-\lambda\right)\left(\sum_{j=1}^m a_j\left(p_j-1\right)^2 t^{p_j-1}\right)}{\left(\sum_{j=1}^m a_j\left(p_j-1\right) t^{p_j-1}\right)^2} .
\end{equation}
For every
$\zeta \in \mathbb{R}$, by direct expansion, we have
\begin{align*}
& \left(\sum_{j=1}^m a_j\left(p_j-1\right) t^{p_j-1}\right)^2 \\
& \quad-\left(\sum_{j=1}^m a_j t^{p_j-1}-\lambda\right)\left(\sum_{j=1}^m a_j\left(p_j-1\right)^2 t^{p_j-1}\right) \\
& =\left(\sum_{j=1}^m a_j\left(p_j-1\right)^2 t^{p_j-1}\right)\left(\lambda-\sum_{j=1}^m a_j\left[1-\zeta\left(p_j-1\right)\right]^2 t^{p_j-1}\right) \\
& \quad+\left(\sum_{j=1}^m a_j\left(p_j-1\right)\left[1-\zeta\left(p_j-1\right)\right] t^{p_j-1}\right)^2 .
\end{align*}
By \eqref{eqkzeta},
$$
\sum_{j=1}^m a_j\left[1-\zeta\left(p_j-1\right)\right]^2 t^{p_j-1} \leq \sum_{j=1}^m a_j k_j t^{p_j-1}<\lambda .
$$
Consequently, the first term on the right-hand side of the identity is strictly positive, since
$$
\sum_{j=1}^m a_j\left(p_j-1\right)^2 t^{p_j-1}>0
$$
while the second term is nonnegative. This proves
$
H_{f, \lambda}^{\prime}(t)>0 .
$

\end{proof}

\begin{lemma}\label{end point}  Let $1<p_1<\cdots<p_m$, and let
$
\phi:[0, \infty) \longrightarrow(0,1]
$
be continuously differentiable, nonincreasing, and convex. Suppose that
$$
p_m-p_1 \leq\left(p_1-1\right) \phi\left(p_m-1\right)+\left(p_m-1\right) \phi\left(p_1-1\right) .
$$
Then there exists $\zeta \in \mathbb{R}$ such that
$$
\left[1-\zeta\left(p_j-1\right)\right]^2 \leq \phi^2 \left(p_j-1\right), ~~ j=1, \ldots, m.
$$
In particular, if $f(t)=\sum_{j=1}^m a_j t^{p_j}$ with $a_j>0$, and
$
\sum_{j=1}^m a_j \phi^2\left(p_j-1\right) t^{p_j-1}<\lambda
$
at some $t>0$, then
$H_{f, \lambda}^{\prime}(t)>0 .$
\end{lemma}

\begin{proof} For each $j$, the inequality
$
\left|1-\zeta\left(p_j-1\right)\right| \leq \phi\left(p_j-1\right)
$
is equivalent to
$$
\frac{1-\phi\left(p_j-1\right)}{p_j-1} \leq \zeta \leq \frac{1+\phi\left(p_j-1\right)}{p_j-1} .
$$
For $q>0$, we have
$1 \ge \phi(0) \geq \phi(q)-q \phi^{\prime}(q)$ by convexity. Hence,
$$
\left(\frac{1-\phi(q)}{q}\right)'=\frac{\phi(q)-q \phi^{\prime}(q)-1}{q^2} \leq0.
$$
Similarly,
$\left(\frac{1+\phi(q)}{q}\right)'<0.$
Thus, for every $j$, we get

$$
\frac{1-\phi\left(p_j-1\right)}{p_j-1} \leq \frac{1-\phi\left(p_1-1\right)}{p_1-1}, \text{and}\quad
\frac{1+\phi\left(p_m-1\right)}{p_m-1} \leq \frac{1+\phi\left(p_j-1\right)}{p_j-1}.
$$
Therefore, it suffices to verify that
$$
\frac{1-\phi\left(p_1-1\right)}{p_1-1} \leq \frac{1+\phi\left(p_m-1\right)}{p_m-1} .
$$
It is equivalent to prove
$p_m-p_1 \leq\left(p_1-1\right) \phi\left(p_m-1\right)+\left(p_m-1\right) \phi\left(p_1-1\right) .$
We may consequently choose any

$$
\zeta \in\left[\frac{1-\phi\left(p_1-1\right)}{p_1-1}, \frac{1+\phi\left(p_m-1\right)}{p_m-1}\right] .
$$
Finally, taking $k_j=\phi^2\left(p_j-1\right)$ in Lemma \ref{sufficient condition}, we  have that  $ H_{f, \lambda}^{\prime}(t)>0$.

\end{proof}

\begin{proposition}\label{prop:peak}
For every compact $J=[a,b]\Subset(1,p_s)$, there exists $C_J>1$, depending
only on $N,s,J$, such that every positive  solution of
\[
(-\Delta)^s v+v=\sum_{j=1}^{m}\theta_jv^{p_j}
\quad\text{in }\R^N,
\qquad
p_j\in J,\quad\theta_j>0,\quad\sum_{j=1}^{m}\theta_j=1,
\]
satisfies $\|v\|_{L^\infty} \le C_J$.
\end{proposition}

\begin{proof}
Suppose by contradiction that there exists a sequence of normalized nonlinearities
$f_n(t)=\sum_{j=1}^{m_n} \theta_{j, n} t^{p_{j, n}}$
with corresponding positive solutions $v_n$ such that
$
M_n:=\|v_n\|_{L^\infty} \to \infty .
$
Without loss of generality, we may assume that $v_n$ is radially symmetric and decreasing.
Define
$\kappa_n=f_n(M_n)/M_n\ge M_n^{a-1} \to +\infty$  and set
\[
 w_n(x)=M_n^{-1}v_n(\kappa_n^{-1/(2s)}x).
\]
Then we have $0<w_n\le1$, $w_n(0)=1$, and
\[
 (-\Delta)^sw_n+\kappa_n^{-1}w_n=\widetilde f_n(w_n),\qquad
 \widetilde f_n(t)=\frac{f_n(M_nt)}{f_n(M_n)}
                =\sum_{j=1}^{m_n} \frac{\theta_{j, n} M_n^{p_{j, n}}}{f_n\left(M_n\right)} t^{p_{j, n}}.
\]
Notice that $\tilde{f}_n$ is again a normalized finite sum of powers, with all exponents belonging to the same fixed interval $[a, b]$. In particular,
$
\tilde{f}_n(0)=0, \quad \tilde{f}_n(1)=1 .
$
Since $a \leq p_{j, n} \leq b$, we have
$t^b \leq \tilde{f}_n(t) \leq t^a$ for $0<t \leq 1,$
and
$t^a \leq \tilde{f}_n(t) \leq t^b$ for  $t \geq 1.$

Computing directly, we have
$$
\tilde{f}_n^{\prime}(t)=\sum_{j=1}^{m_n} \frac{\theta_{j, n} M_n^{p_{j, n}}}{f_n\left(M_n\right)} p_{j, n} t^{p_{j, n}-1} .
$$
Consequently, for every $T>0$,
$0 \leq \tilde{f}_n^{\prime}(t) \leq b \max \left\{1, T^{b-1}\right\}$ for $0 \leq t \leq T$.
Hence the sequence $\{\tilde{f}_n\}$ is uniformly bounded and uniformly Lipschitz on every bounded interval. By the Arzelà-Ascoli theorem and a diagonal subsequence argument, there exists a continuous function $f_{\infty}$ such that
$$
\tilde{f}_n \to f_{\infty} \quad \text { uniformly on every bounded interval of }[0, \infty) \text {. }
$$
By interior regularity estimates for the fractional Laplacian,  we have that the sequence \(w_n\) is locally equi-continuous.
After passing to a subsequence,
$w_n \to w $ locally uniformly in $\mathbb{R}^N$.
In particular,

$$
(-\Delta)^s w=f_{\infty}(w) \quad \text { in } \mathbb{R}^N, \quad 0 \leq w \leq 1, \quad\text{and}\quad w(0)=1 .
$$

When $N>2s$, we notice that
\[
 \frac{tf_\infty'(t)}{f_\infty(t)}\in[a,b],\qquad
 \frac{t^2f_\infty''(t)}{f_\infty(t)}
 \le b(b-1).
\]
Thus $f_\infty$ is increasing and  $f'_\infty$ is locally Lipschitz, and the regularity argument yields $w \in C_{\text {loc }}^{1,1}$. Moreover
\[
 \int_{\R^N}\frac{|w(x)|}{1+|x|^{N+2s}}\,dx<\infty,
\]
which follows directly from $w\le1$.
Since
\[
 \frac{d}{dt}\left(\frac{f_\infty(t)}{t^{p_s}}\right)
 =\frac{f_\infty(t)}{t^{p_s+1}}
   \left(\frac{tf_\infty'(t)}{f_\infty(t)}-p_s\right)\leq\frac{f_\infty(t)}{t^{p_s+1}}
   \left(b-p_s\right)<0,
\]
we obtain that $f_\infty(t)/t^{p_s}$ is strictly decreasing.
Then the moving-spheres
classification \cite[Theorem 3.1]{CLZ} excludes this bounded positive
solution.

When $N\le2s$, it holds that $N=1$, $s\ge1/2$. Since $w(0)=1$ and $w$ is continuous, there is a ball $B$ centered at zero such that for all $ x \in B$,
$w(x) \geq \frac{1}{2} .$
Since $f_{\infty}$ is increasing, for every $x \in B$ we have
$f_{\infty}(w(x)) \geq f_{\infty}(1 / 2) \geq 2^{-b} .$
Thus, with $d=2^{-b}>0$,
$f_{\infty}(w) \geq d  \chi_B$ in $\mathbb{R}^N$, where $\chi_B$ denotes the characteristic function of $B$.

For $\varepsilon>0$, let $G_{\varepsilon}$ be the positive kernel of
$
\left((-\Delta)^s+\varepsilon\right)^{-1}.
$
Adding $\varepsilon w$ to the limiting equation, we have $$\left((-\Delta)^s+\varepsilon\right) w=f_{\infty}(w)+\varepsilon w \geq d \chi_B .$$
Therefore 
$w(x) \geq d \int_B G_{\varepsilon}(x-y) d y .$
In particular, we have
\begin{equation}\label{wo}
1=w(0) \geq d \int_B G_{\varepsilon}(y) d y .
\end{equation}
We will show that the right-hand side  of (\ref{wo}) tends to infinity as $\varepsilon \to 0$.

Recall that
$$
G_{\varepsilon}(y)=\int_0^{\infty} e^{-\varepsilon \tau} P_s(\tau, y) d \tau,
$$
where the fractional heat kernel satisfies
$
P_s(\tau, y)=\tau^{-N /(2 s)} P_s\left(1, \tau^{-1 /(2 s)} y\right) .
$
For $y$ in the fixed bounded ball $B$,
$
\tau^{-1 /(2 s)} y \longrightarrow 0
$
uniformly as $\tau \rightarrow \infty$. Since $P_s(1, \cdot)$ is continuous and positive at zero, there are constants $c>0$ and $T>0$ such that

$$
P_s(\tau, y) \geq c \tau^{-N /(2 s)} \quad(\forall y \in B, \tau \geq T) .
$$
For sufficiently small $\varepsilon$, we have $1 / \varepsilon>T$. On the interval $T \leq \tau \leq 1 / \varepsilon$,
$
e^{-\varepsilon \tau} \geq e^{-1} .
$
Hence
$$
\begin{aligned}
1 \geq d \int_B G_{\varepsilon}(y) d y 
 \geq d c|B| e^{-1} \int_T^{1 / \varepsilon} \tau^{-N /(2 s)} d \tau.
\end{aligned}
$$
If $N=2 s$, then

$$
\int_T^{1 / \varepsilon} \frac{d \tau}{\tau}=\log \frac{1}{\varepsilon T} \longrightarrow \infty .
$$
If $N<2 s$, then
$$
\int_T^{1 / \varepsilon} \tau^{-N /(2 s)} d \tau=\frac{\varepsilon^{-(1-N /(2 s))}-T^{1-N /(2 s)}}{1-N /(2 s)} \longrightarrow \infty .
$$
Both cases contradict $w(0)=1$. Thus the limiting solution cannot exist when $N \leq 2 s$. This completes the proof.
\end{proof}

For a smaller range of exponents, the constant in
Proposition~\ref{prop:peak} admits the following explicit bound.

\begin{proposition}\label{prop:peak-explicit}
Let $N\geq 2$, $0<s<1$, and $J=[a,b]$ with
$1<a\leq b<N/(N-2s).$
Let $v\in H^s(\mathbb R^N)\cap L^\infty(\mathbb R^N)$ be a positive
weak solution of
\begin{equation}\label{eq:peak-normalized}
\begin{gathered}
(-\Delta)^s v+v=\hat f(v)\quad\text{in }\mathbb R^N,\\
\end{gathered}
\end{equation}
where $ \hat f(t)=\sum_{j=1}^{m}\theta_j t^{p_j}, 
 p_j\in J, \theta_j>0,$ and $ \sum_{j=1}^{m}\theta_j=1.$
Set $k=N/(2s)$. 
Then, for every $T\geq\log k$,
\begin{equation}\label{eq:peak-explicit-bound}
\|v\|_{L^\infty}\leq C_J(T)
:=
\left[
\frac{kT}{k-1}
\left(1+\frac{k}{T}\right)^{\frac{a}{k(a-1)}}
\right]^{\frac{k}{k-b(k-1)}}.
\end{equation}
The minimum of $C_J(T)$ over $T\geq\log k$ is attained at
\[
T_{\log}=\max\left\{\log k,\frac{a}{a-1}-k\right\}.
\]
In particular, we may take $C_J=C_J(T_{\log})$.

\end{proposition}

\begin{proof}
Let $P_s(\tau,x)$ be the heat kernel of $(-\Delta)^s$. Then $P_s$ is positive, has integral one, and
is radially nonincreasing in $x$. The scaling gives
\begin{equation}\label{eq:peak-kernel}
\tau^k P_s(\tau,x)
=P_s\bigl(1,\tau^{-1/(2s)}x\bigr)
\leq P_s\bigl(1,T^{-1/(2s)}x\bigr)
=T^k P_s(T,x)
\end{equation}
for $0<\tau\leq T,  x\in\mathbb R^N$.

Set $M=\|v\|_{L^\infty}$. Since
$\hat f(v)\leq (\sum_j\theta_j M^{p_j-1})v$, we have
$\hat f(v)\in L^2(\mathbb R^N)$. For $\tau>0$, define
\[
Y_\tau=P_s(\tau,\cdot)*v,
\quad \text{and} \quad
X_\tau=P_s(\tau,\cdot)*\hat f(v).
\]
Using
$\partial_\tau Y_\tau=-(-\Delta)^s Y_\tau$
and \eqref{eq:peak-normalized},
we have $X_T=Y_T-\partial_TY_T$ in $L^2(\mathbb R^N)$.
By \eqref{eq:peak-kernel} and $v\geq 0$, the function
$T^kY_T$ is nondecreasing in $T$. Hence, we have
\begin{equation}\label{eq:peak-semigroup}
\partial_\tau Y_\tau\geq -\frac{k}{\tau}Y_\tau,
\quad \text{and}\quad
X_T\leq \left(1+\frac{k}{T}\right)Y_T.
\end{equation}
Since $\hat f$ is convex and
$P_s(T,\cdot)\,$ has integral one, by 
Jensen's inequality, we have
\[
\begin{aligned}
\hat f(Y_T(x))
&= \hat f\!\left(\int_{\mathbb R^N}P_s(T,x-y)v(y)\,dy\right)\\
&\leq \int_{\mathbb R^N}P_s(T,x-y)\hat f (v(y))\,dy
= X_T(x).
\end{aligned}
\]
When $Y_T\geq 1$, we have $\hat f(Y_T)\geq Y_T^a$ because
$p_j\geq a$ and $\sum_j\theta_j=1$. Combining this with
\eqref{eq:peak-semigroup}, we have
\[
Y_T\leq \left(1+\frac{k}{T}\right)^{\frac{1}{a-1}}.
\]
The same bound holds when $Y_T<1$, since $\left(1+\frac{k}{T}\right)^{\frac{1}{a-1}}$ is greater
than one. Therefore, we get
\begin{equation}\label{eq:peak-average}
\|X_T\|_{L^\infty}\leq \left(1+\frac{k}{T}\right)^{\frac{a}{a-1}}.
\end{equation}

Notice that
\begin{equation}\label{eq:peak-resolvent}
v(x)=G_{s,1}*\hat f(v)=\int_0^\infty e^{-\tau}X_\tau(x)\,d\tau
\quad\text{for almost every }x\in\mathbb R^N.
\end{equation}
Since $G_{s,1}$ has integral one, we have $M\leq \hat f(M)$. If $0<M<1$,
then $\hat f (M)\leq M^a<M$, a contradiction. Therefore, $M\geq 1$ and
$\hat f (v)\leq M^b$. By
\eqref{eq:peak-kernel} and \eqref{eq:peak-average},
we have
\[
0\leq X_\tau\leq \min\left\{M^b,\left(1+\frac{k}{T}\right)^{\frac{a}{a-1}} (T/\tau)^k\right\}
\qquad\text{for }0<\tau\leq T.
\]
For $\tau>0$, set 
\[
\Phi_T(\tau)=\min\left\{M^b,\left(1+\frac{k}{T}\right)^{\frac{a}{a-1}}(T/\tau)^k\right\}.
\]
Therefore, $X_\tau\leq\Phi_T(\tau)$ for $0<\tau\leq T$.

For $\tau>T$, 
$X_\tau=P_s(\tau-T,\cdot)*X_T$. Since $X_T\leq\Phi_T(T)$, we have
\[
X_\tau(x)
=\int_{\mathbb R^N}P_s(\tau-T,x-y)X_T(y)\,dy
\leq \Phi_T(T).
\]
Using \eqref{eq:peak-resolvent}
and $e^{-\tau}\leq 1$ on $(0,T]$, we obtain
\begin{equation}\label{eq:peak-tail}
\begin{aligned}
M
&\leq \int_0^T\Phi_T(\tau)\,d\tau
+\Phi_T(T)\int_T^\infty e^{-\tau}\,d\tau\\
&= \int_0^T\Phi_T(\tau)\,d\tau+e^{-T}\Phi_T(T).
\end{aligned}
\end{equation}

For $\tau\geq T$,  $(T/\tau)^k\leq 1$. Hence
\[
\begin{aligned}
\Phi_T(\tau)
&=\min\left\{M^b, \left(1+\frac{k}{T}\right)^{\frac{a}{a-1}}\left(\frac{T}{\tau}\right)^k\right\}\\
&\geq \min\left\{M^b\left(\frac{T}{\tau}\right)^k,
                 \left(1+\frac{k}{T}\right)^{\frac{a}{a-1}}\left(\frac{T}{\tau}\right)^k\right\}\\
&=\left(\frac{T}{\tau}\right)^k\Phi_T(T).
\end{aligned}
\]
Since $k>1$, we have
\begin{equation}\label{eq:peak-tail-lower}
\int_T^\infty\Phi_T(\tau)\,d\tau
\geq \Phi_T(T)\int_T^\infty\left(\frac{T}{\tau}\right)^k\,d\tau
=\frac{T}{k-1}\Phi_T(T).
\end{equation}
Assume $T\geq\log k$, then we have
\[
\log k=\int_1^k\frac{dt}{t}\geq \frac{k-1}{k},
\quad \text{and}\quad
e^{-T}\leq \frac{1}{k}\leq \frac{\log k}{k-1}
\leq \frac{T}{k-1}.
\]
Therefore $e^{-T}\Phi_T(T)\leq\int_T^\infty\Phi_T(\tau)\,d\tau$ by
\eqref{eq:peak-tail-lower}. Consequently, by
\eqref{eq:peak-tail} we have 
\begin{equation}\label{eq:peak-integral}
M\leq \int_0^\infty\Phi_T(\tau)\,d\tau.
\end{equation}

Notice that  $M^b$ and $\left(1+\frac{k}{T}\right)^{\frac{a}{a-1}}(T/\tau)^k $ agree at
\[
\tau_0=T\left(\frac{\left(1+\frac{k}{T}\right)^{\frac{a}{a-1}}}{M^b}\right)^{1/k}.
\]
Therefore, 
\[
\begin{aligned}
\int_0^\infty\Phi_T(\tau)\,d\tau
&=M^b\tau_0+\frac{\left(1+\frac{k}{T}\right)^{\frac{a}{a-1}}T^k}{k-1}\tau_0^{1-k}\\
&=\frac{k}{k-1}M^b\tau_0\\
&=\frac{kT}{k-1}\left(1+\frac{k}{T}\right)^{\frac{a}{k(a-1)}} M^{b(k-1)/k}.
\end{aligned}
\]

By \eqref{eq:peak-integral}, we have
\[
M^{1-b(k-1)/k}\leq \frac{kT}{k-1}\left(1+\frac{k}{T}\right)^{\frac{a}{k(a-1)}}.
\]
Since $b<N/(N-2s)=k/(k-1)$, we obtain $M\leq C_J(T)$.

Finally, differentiating \eqref{eq:peak-explicit-bound}, we have
\[
\frac{d}{dT}\log C_J(T)
=\frac{k}{k-b(k-1)}
\frac{T+k-a/(a-1)}{T(T+k)}.
\]
The derivative is negative for $T<a/(a-1)-k$ and positive for
$T>a/(a-1)-k$. If the critical point lies below $\log k$, the minimum
on $[\log k,\infty)$ is attained at $\log k$. Otherwise, it is attained
at $T=a/(a-1)-k$. This gives
$T_{\log}=\max\{\log k,a/(a-1)-k\}$.

\end{proof}

\begin{lemma}\label{thm:uniform}
Let $N\ge1$ and $0<s<1$. For a compact interval $J\Subset(1,p_s)$, we suppose that
\[
f(t)=\sum_{j=1}^{m}a_jt^{p_j},
\qquad a_j>0,\qquad1<p_1<\cdots<p_m,\qquad p_j\in J,
\]
and
\[
p_m-p_1\le(p_1-1)C_J^{-(p_m-1)/2}
+(p_m-1)C_J^{-(p_1-1)/2},
\]
where $C_J > 1$ is the constant given in Proposition \ref{prop:peak}, depending only on $N,s,J$.
Then, for every  positive solution $u \in H^s$ of
$(-\Delta)^s u+\lambda u=f(u)$ in $\R^N,$
the function
$H_{f,\lambda}(t)$
satisfies
$$
H_{f,\lambda}'(t)>0, \quad\text{for}\quad 0<t<\|u\|_{L^\infty}.
$$
Moreover, if $\rho_\lambda>0$ denotes the unique positive solution of
$
\sum_{j=1}^{m}a_j\rho_\lambda^{p_j-1}=\lambda,
$
then
\[
\|u\|_{L^\infty}\le C_J\rho_\lambda.
\]
\end{lemma}

\begin{proof}
The function $f(t)/t=\sum_ja_jt^{p_j-1}$ is continuous and strictly
increasing from zero to infinity. Thus there is a unique
$\rho_\lambda>0$ satisfying $f(\rho_\lambda)=\lambda\rho_\lambda$.
Set $v(y)=\rho_\lambda^{-1}u(\lambda^{-1/(2s)}y)$ and
$\theta_j=a_j\rho_\lambda^{p_j-1}/\lambda$. Then $\sum_j\theta_j=1$ and
$(-\Delta)^sv+v=\sum_j\theta_jv^{p_j}$. Proposition~\ref{prop:peak}
gives $\|u\|_{L^\infty}\le C_J\rho_\lambda$. 

For $0<t<\|u\|_{L^\infty}$, take $k_j=C_J^{-(p_j-1)}$. Thus
\begin{equation}\label{stri}
 \sum_j  a_j k_j t^{p_j-1}/\lambda=\sum_j\theta_j
       \left(\frac{t}{C_J\rho_\lambda}\right)^{(p_j-1)}<1.
\end{equation}
Take
$\phi(q)=C_J^{-q/2}, C_J>1$ .
This function is decreasing and convex, with $ \phi(0)=1$. 
Then, Lemma \ref{end point} yields
$$
H_{f, \lambda}^{\prime}(t)>0 \quad\text{for}\quad 0<t<\|u\|_{L^\infty}.
$$

\end{proof}

\subsection{The case of sums of powers: $N = 1$}

When $N=1$, we provide different conditions.

\begin{lemma}\label{A}
Under the assumptions 
$1<p_i, p_j \leq \frac{4}{C}-1$ for some $1\leq C<2 ,$ there holds
\begin{align}\label{inequ}
(2-C)\left(p_i-p_j\right)^2 \leq & \left(\frac{4}{p_i+1}-C\right)\left(p_j-1\right)^2 
+\left(\frac{4}{p_j+1}-C\right)\left(p_i-1\right)^2 .
\end{align}

\end{lemma}

\begin{proof}
Since the inequality is symmetric in $p_i$ and $p_j$, we may assume that $p_i \geq p_j$. Using

$$
\left(p_i-1\right)^2+\left(p_j-1\right)^2=\left(p_i-p_j\right)^2+2\left(p_i-1\right)\left(p_j-1\right),
$$
we obtain

$$
\begin{aligned}
& \left(\frac{4}{p_i+1}-C\right)\left(p_j-1\right)^2+\left(\frac{4}{p_j+1}-C\right)\left(p_i-1\right)^2-(2-C)\left(p_i-p_j\right)^2 \\
& =2\left(p_i-1\right)\left(p_j-1\right)\left(\frac{4}{p_i+1}-C\right) \\
& \quad+\left(\frac{4}{p_j+1}-\frac{4}{p_i+1}\right)\left(p_i-1\right)^2+\left(\frac{4}{p_i+1}-2\right)\left(p_i-p_j\right)^2 .
\end{aligned}
$$
 Indeed,
$$
\begin{aligned}
& \left(\frac{4}{p_j+1}-\frac{4}{p_i+1}\right)\left(p_i-1\right)^2+\left(\frac{4}{p_i+1}-2\right)\left(p_i-p_j\right)^2 \\
& =\frac{2\left(p_i-1\right)\left(p_j-1\right)\left(p_i-p_j\right)\left(p_j+2-p_i\right)}{\left(p_i+1\right)\left(p_j+1\right)}.
\end{aligned}
$$
Therefore,
$$
\begin{aligned}
& \left(\frac{4}{p_i+1}-C\right)\left(p_j-1\right)^2+\left(\frac{4}{p_j+1}-C\right)\left(p_i-1\right)^2-(2-C)\left(p_i-p_j\right)^2 \\
& =2\left(p_i-1\right)\left(p_j-1\right)\left[\frac{4}{p_i+1}-C+\frac{\left(p_i-p_j\right)\left(p_j+2-p_i\right)}{\left(p_i+1\right)\left(p_j+1\right)}\right] .
\end{aligned}
$$

Both terms inside the brackets are nonnegative. First, the assumption $p_i \leq 4 / C-1$ yields $\frac{4}{p_i+1}-C \geq 0 .$
Then, since $C \geq 1$, $p_i \leq \frac{4}{C}-1 \leq 3,$ and hence $p_j+2-p_i >0 .$

Together with $p_i-p_j \geq 0$, the desired inequality follows. Moreover, it is strict whenever $p_i \neq p_j$. 
\end{proof}

\begin{lemma}\label{thm:line}
Let $N=1$, $f(t)=\sum_{j=1}^m a_jt^{p_j}$, and let $u$ be a positive,
even, strictly decreasing solution of \eqref{Q1}.
Under either condition $(i)$ or $(ii)$ of Theorem~\ref{thmN1},
$H_{f,\lambda}'(t)>0$ for $0<t\le u(0)$.
Under condition $(iii)$, $H_{f,\lambda}'(t)>0$ for $0<t<u(0)$.
\end{lemma}

\begin{proof}

Put $M=u(0)$, $g(t)=f(t)-\lambda t$, and
$\mathcal G(t)=\int_0^tg(\tau)\,d\tau=F(t)-\lambda t^2/2$.
Define $U(x,y)=(P_{s,y}*u)(x), y>0$,
where
\[
 P_{s,y}(x)=\frac{\Gamma(s+1/2)}{\sqrt\pi\Gamma(s)}
 \frac{y^{2s}}{(x^2+y^2)^{(1+2s)/2}}.
\]
The extension satisfies $\operatorname{div}\left(y^{1-2 s} \nabla U\right)=0, \quad U(x, 0)=u(x),$
and
\begin{equation}\label{V0y}
-d_s^{\text {ext }} \lim _{y \downarrow 0} y^{1-2 s} U_y(x, y)=f(u(x))-\lambda u(x), \qquad d_s^{\mathrm{ext}}=\frac{2^{2 s-1} \Gamma(s)}{\Gamma(1-s)}>0 .
\end{equation}
Since $u$ is even, $U$ is also even in $x$, so $U_x(0, y)=0$.

Consider
\begin{equation}\label{Ha}
\frac{d_s^{\text {ext }}}{2} \int_0^{\infty} y^{1-2 s}\left(U_x^2-U_y^2\right)(x, y) d y+F(u(x))-\frac{\lambda}{2} u(x)^2.
\end{equation}
We claim that the  expression (\ref{Ha}) is independent of $x$. Indeed,
\begin{align}\label{d=0}
\frac{d}{d x} & {\left[\frac{d_s^{\text {ext }}}{2} \int_0^{\infty} y^{1-2 s}\left(U_x^2-U_y^2\right) d y\right] } 
=d_s^{\text {ext }} \int_0^{\infty} y^{1-2 s}\left(U_x U_{x x}-U_y U_{x y}\right) d y.
\end{align}
The extension equation gives
$y^{1-2 s} U_{x x}=-\partial_y\left(y^{1-2 s} U_y\right),$ and
integrating by parts, we have
$$
-d_s^{\operatorname{ext}}\left[y^{1-2 s} U_x U_y\right]_{y=0}^{y=\infty}=-[f(u)-\lambda u] u^{\prime} .
$$
On the other hand,
$$
\frac{d}{d x}\left(F(u)-\frac{\lambda}{2} u^2\right)=[f(u)-\lambda u] u^{\prime},
$$
therefore,  we have $\frac{d}{d x}\left[\frac{d_s^{\mathrm{ext}}}{2} \int_0^{\infty} y^{1-2 s}\left(U_x^2-U_y^2\right) d y+F(u)-\frac{\lambda}{2} u^2\right]=0$.

Evaluating  (\ref{Ha})  at $x=0$, where $U_x(0, y)=0$, we have
\begin{equation}\label{FM}
F(M)-\frac{\lambda M^2}{2}=\frac{d_s^{\mathrm{ext}}}{2} \int_0^{\infty} y^{1-2 s} U_y(0, y)^2 d y.
\end{equation}

Set $V(x, y):=-y^{1-2 s} U_y(x, y).$
The extension equation implies
$V_y=y^{1-2 s} U_{x x},$ and $ V_x=-y^{1-2 s} U_{x y}.$
Consequently, we have $\operatorname{div}\left(y^{2 s-1} \nabla V\right)=0.$
Since $2 s-1=1-2(1-s)$, this is precisely the extension equation of order $1-s$. Therefore, by the Poisson representation $V$ is the $(1-s)$-Poisson extension of its boundary trace. By (\ref{V0y}), we have
$V(x, 0)=\frac{f(u(x))-\lambda u(x)}{d_s^{\mathrm{ext}}}$.
Thus,
$$V(0, y)=\frac{1}{d_s^{\mathrm{ext}}} \int_{\mathbb{R}} P_{1-s, y}(x)[f(u(x))-\lambda u(x)] d x.$$

Note that
$
f(M)-\lambda M=(-\Delta)^s u(0)=c_{1, s} \operatorname{P.V.} \int_{\mathbb{R}} \frac{M-u(x)}{|x|^{1+2 s}} d x>0 .
$

On the one hand, we will show that
$U_y(0, y)<0$ and $ V(0, y)>0 .$
Using the scaling of the Poisson kernel
$U(0, y)=\int_{\mathbb{R}} P_{s, 1}(z) u(y z) d z,$
a computation yields
$
U_y(0, y)=\int_{\mathbb{R}} P_{s, 1}(z) z u^{\prime}(y z) d z<0.
$
The positivity of $V(0, y)$ follows from $V(x, y):=-y^{1-2 s} U_y(x, y).$

On the other hand, the function $t \mapsto f(t)-\lambda t$ is strictly convex and vanishes at zero. Hence
$
f(t)-\lambda t<\frac{t}{M}[f(M)-\lambda M]$ for $0<t<M .$
Integrating this inequality against $P_{1-s, y}$, we have
$$
V(0, y)<\frac{f(M)-\lambda M}{M d_s^{\mathrm{ext}}}\left(P_{1-s, y} * u\right)(0) .
$$
We claim that
$\left(P_{1-s, y} * u\right)(0) \leq C U(0, y),$
where $U(0, y)=\left(P_{s, y} * u\right)(0) ,$ and
$C=1$ if $ s \geq \frac{1}{2}$,  $C=C_s= \frac{\Gamma\left(\frac{3}{2}-s\right) \Gamma(s)}{\Gamma(1-s) \Gamma\left(s+\frac{1}{2}\right)}$ if $ s<\frac{1}{2}$. For $s<1 / 2$, this claim follows from 
$$
\frac{P_{1-s, y}(x)}{P_{s, y}(x)}=C_s\left(1+\frac{x^2}{y^2}\right)^{2 s-1} \leq C_s .
$$
For $s>1 / 2$, the two kernels have equal total mass and cross once: $P_{s, y}$ is larger near zero and smaller farther away. Integration against the even decreasing function $u$ therefore gives the comparison with constant one. At $s=1 / 2$, the kernels coincide.

Hence, we obtain

\begin{equation}\label{V}
0<V(0, y)<\frac{C[f(M)-\lambda M]}{M d_s^{\mathrm{ext}}} U(0, y) .
\end{equation}
Since
$y^{1-2 s} U_y^2=\left(-U_y\right) V,$
equation (\ref{FM}) becomes

$$
F(M)-\frac{\lambda M^2}{2}=\frac{d_s^{\mathrm{ext}}}{2} \int_0^{\infty}\left[-U_y(0, y)\right] V(0, y) d y
$$
Using  (\ref{V}) and $-U_y>0$, we have

$$
F(M)-\frac{\lambda M^2}{2}<\frac{C[f(M)-\lambda M]}{2 M} \int_0^{\infty}\left[-U_y(0, y)\right] U(0, y) d y .
$$
Clearly,
$\int_0^{\infty}\left[-U_y(0, y)\right] U(0, y) d y=-\frac{1}{2}\left[U(0, y)^2\right]_0^{\infty}=\frac{M^2}{2},$
since $U(0,0)=M$ and $U(0, y) \to 0$ as $y \to \infty.$
Therefore,
$$4\left(F(M)-\frac{\lambda M^2}{2}\right)<C M[f(M)-\lambda M] .$$
For the power mixture, this is exactly

$$
\sum_{j=1}^m a_j\left(\frac{4}{p_j+1}-C\right) M^{p_j-1}<(2-C) \lambda .
$$
Since
$p_j>1$, for every $0<t \leq M$, we obtain
\begin{equation}\label{2cl}
\sum_{j=1}^m a_j\left(\frac{4}{p_j+1}-C\right) t^{p_j-1}<(2-C) \lambda  .
\end{equation}

Next, we will  show $H_{f, \lambda}^{\prime}(t)>0$ for $0<t\le M$.
A direct calculation yields
$$
t H_{f, \lambda}^{\prime}(t)=\frac{\lambda \sum_{j=1}^m a_j\left(p_j-1\right)^2 t^{p_j-1}-\sum_{i<j} a_i a_j\left(p_i-p_j\right)^2 t^{p_i+p_j-2}}{\left(\sum_{j=1}^m a_j\left(p_j-1\right) t^{p_j-1}\right)^2} .
$$
Since the denominator  is positive, we only need to prove
\begin{equation}\label{coss}
\sum_{i<j} a_i a_j\left(p_i-p_j\right)^2 t^{p_i+p_j-2}<\lambda \sum_j a_j\left(p_j-1\right)^2 t^{p_j-1} . 
\end{equation}

Assume that
\begin{align} \label{eq4/C-1}
	1<p_i, p_j \leq \frac{4}{C}-1.
\end{align}
By Lemma \ref{A}
we have  
\begin{align}\label{inequ:center}
(2-C)\left(p_i-p_j\right)^2 \leq & \left(\frac{4}{p_i+1}-C\right)\left(p_j-1\right)^2 
+\left(\frac{4}{p_j+1}-C\right)\left(p_i-1\right)^2 .
\end{align}
Multiplying (\ref{inequ:center}) by 
$a_i a_j t^{p_i+p_j-2}$ yields
$$
\begin{aligned}
& (2-C) \sum_{i<j} a_i a_j\left(p_i-p_j\right)^2 t^{p_i+p_j-2} \\
& \leq \sum_{i<j} a_i a_j t^{p_i+p_j-2}\left[\left(\frac{4}{p_i+1}-C\right)\left(p_j-1\right)^2+\left(\frac{4}{p_j+1}-C\right)\left(p_i-1\right)^2\right]\\
& \leq\left[\sum_{i=1}^m a_i\left(\frac{4}{p_i+1}-C\right) t^{p_i-1}\right]\left[\sum_{j=1}^m a_j\left(p_j-1\right)^2 t^{p_j-1}\right] \\
& <(2-C) \lambda \sum_{j=1}^m a_j\left(p_j-1\right)^2 t^{p_j-1},
\end{aligned}
$$
where the strict inequality follows from (\ref{2cl}) and the positivity of the second factor. Since $2-C>0$, (\ref{coss}) holds.

For condition $(i)$, $s \geq \frac{1}{2}$, \eqref{eq4/C-1} indeed holds with $C=1$; while for condition $(ii)$, $0<s<\frac{1}{2}$, \eqref{eq4/C-1} indeed holds with $C=C_s$.
We need $C_s<2$.
To determine when this holds, we observe that
$$
\frac{d}{d s} \log C_s=\psi(s)-\psi\left(s+\frac{1}{2}\right)+\psi(1-s)-\psi\left(\frac{3}{2}-s\right),
$$
where $\psi=\Gamma^{\prime} / \Gamma$. The digamma function $\psi$ is strictly increasing on  $(0, \infty)$, implying
$C_s^{\prime}<0 .$
Moreover,
$\lim _{s \to 0} C_s=+\infty, \text{and} \lim _{s \to 1 / 2} C_s=1 .$
Hence, there is exactly one $s_0 \in\left(0, \frac{1}{2}\right)$ such that
$C_{s_0}=2 .$
For $s_0<s<\frac{1}{2}$, we have $1<C_s<2$.

\vskip0.15in
In the following we prove the result under condition $(iii)$. 
For $s>1/2$, set $h_s=G_s/c_s$. By the Stieltjes representation of
$(1+z^s)^{-1}$, we have
\[
 G_s(x)=\frac{\sin(\pi s)}\pi\int_0^\infty
 \frac{r^{2s}e^{-r|x|}}{1+2r^{2s}\cos(\pi s)+r^{4s}}\,dr,
\]
and 
$$
c_s=G_s(0)=\frac{1}{\pi} \int_0^{\infty} \frac{d \xi}{1+\xi^{2 s}}=\frac{1}{2 s \sin (\pi /(2 s))} .
$$
Note that $h_s(0)=1$,
$0<h_s(x) \leq 1$ and $h_s$ is strictly decreasing in $x > 0$.

Observe that
\begin{align*}
& h_s(A+B)-h_s(A) h_s(B) \\
& =\frac{1}{2}\left(\frac{\sin (\pi s)}{\pi c_s}\right)^2 \int_0^{\infty} \int_0^{\infty} \frac{r^{2 s} q^{2 s}\left(e^{-A r}-e^{-A q}\right)\left(e^{-B r}-e^{-B q}\right)}{\left[1+2 r^{2 s} \cos (\pi s)+r^{4 s}\right]\left[1+2 q^{2 s} \cos (\pi s)+q^{4 s}\right]} d r d q.
\end{align*}
It is clear that
$\left(e^{-A r}-e^{-A q}\right)\left(e^{-B r}-e^{-B q}\right) \geq 0.$
Hence, we have
$h_s(A+B) \geq h_s(A) h_s(B)$ for $A, B \geq 0.$ Since $|x-y| \leq|x|+|y|$, monotonicity of $h_s$ then gives
$h_s(x-y) \geq h_s(|x|+|y|) \geq h_s(x) h_s(y) .$
Equivalently,
\begin{equation}\label{nor}
G_s(x-y) \geq G_s(x) G_s(y)/c_s .
\end{equation}

Set $w(x)=u\left(\lambda^{-1 /(2 s)} x\right) .$
Then $w(0)=M$, and  $(-\Delta)^s w+w=f(w)/\lambda.$
Hence
\begin{equation}\label{wx}
w(x)=\frac{1}{\lambda} \int_{\mathbb{R}} G_s(x-y) f(w(y)) d y .
\end{equation}
Using (\ref{nor}), we have
\begin{equation}\label{w2}
w(x) \geq \frac{G_s(x)}{\lambda c_s} \int_{\mathbb{R}} G_s(y) f(w(y)) d y .
\end{equation}
Evaluating (\ref{wx}) at zero, we obtain
$$
M=\frac{1}{\lambda} \int_{\mathbb{R}} G_s(x) f(w(x)) d x.
$$
Therefore, by (\ref{w2}), we have
\begin{equation}\label{w3}
w(x) \geq M h_s(x) .
\end{equation}
For the positive power mixture, $f$ is increasing. By (\ref{w3}),

$$
\begin{aligned}
M & \geq \frac{1}{\lambda} \int_{\mathbb{R}} G_s(x) f\left(M h_s(x)\right) d x =\frac{1}{\lambda} \sum_{j=1}^m a_j M^{p_j} \int_{\mathbb{R}} G_s(x) h_s^{p_j}(x) d x.
\end{aligned}
$$
Define 
$K_{s, r}:=\int_{\mathbb{R}} G_s(x) h_s^r (x)d x.$
Dividing by $M>0$, we obtain
\begin{equation}\label{sum}
\sum_{j=1}^m a_j K_{s, p_j} M^{p_j-1} \leq \lambda .
\end{equation}
Since $p_j>1$, it follows that
$\sum_{j=1}^m a_j K_{s, p_j} t^{p_j-1}<\lambda$ for $ 0<t<M .$

Take $k_j=K_{s, p_j}$ as in Lemma \ref{sufficient condition}.
We want to apply Lemma \ref{end point} with $
\phi(z)=\sqrt{K_{s, 1+z}}$ for $z \geq 0.$
Recall that
$\int_{\mathbb{R}} G_s(x) d x=1$ and $0<h_s \leq 1$.
Thus $0<K_{s, r} \leq 1.$
By the monotonicity of $h_{s}^r$, we get that $K_{s, r}$ is nonincreasing.
By Hölder's inequality
$$
K_{s, \theta r+(1-\theta) q} \leq K_{s, r}^\theta K_{s, q}^{1-\theta} \quad\text{for }~ 0 \leq \theta \leq 1 .
$$
Taking square roots, we have
$$
\begin{aligned}
\sqrt{K_{s, \theta r+(1-\theta) q}} & \leq\left(\sqrt{K_{s, r}}\right)^\theta\left(\sqrt{K_{s, q}}\right)^{1-\theta}  \leq \theta \sqrt{K_{s, r}}+(1-\theta) \sqrt{K_{s, q}} .
\end{aligned}
$$
Thus $\phi$ is convex. Its continuous differentiability follows by differentiation under the integral, using
$h_s^{1+z}\left|\log h_s\right| \leq \sup _{0<c \leq 1} c|\log c|<\infty$
and integrability of $G_s$.

Lemma \ref{end point} therefore applies whenever
$$
p_m-p_1 \leq\left(p_1-1\right) \sqrt{K_{s, p_m}}+\left(p_m-1\right) \sqrt{K_{s, p_1}},
$$
and we derive
$H_{f, \lambda}^{\prime}(t)>0$ for $0<t<M.$
\end{proof}

\subsection{The case involving a logarithmic term} 

We first prove a result with general nonlinearities.
 
\begin{proposition}
\label{thm:general-quotient-monotonicity}
Let $N\geq 1$,  $0<s<1$ and $\lambda>0$.
Suppose that
\[
f\in C([0,\infty))\cap C^2((0,\infty)),
\quad f(0)=0,\quad f(t)>0~~ \text{for}~~ t>0.
\]
Define
\[
\eta(t):=\frac{t f'(t)}{f(t)}-1,
\qquad D:=t\frac{d}{dt},
\]
and assume that, for every $t>0$,
\[
\eta(t)>0,
\qquad D^2\!\left(\frac{f(t)}{t}\right)\geq 0.
\]
Let $u\in H^s(\mathbb R^N)\cap L^\infty(\mathbb R^N)$ be a positive
 solution of
$(-\Delta)^s u+\lambda u=f(u)$ in $\mathbb R^N,$
and set $M:=\|u\|_{L^\infty}$.
We further assume that there exists a constant $\hat C>1$ such that
$\frac{f(M)}{\lambda M}\leq \hat C$ and $\chi(\hat C-1)\leq 1$, where
\[
\chi:=\sup_{t>0}\frac{(D\eta(t))_+}{\eta(t)^2}.
\]

Then $H_{f, \lambda}'(t)>0$ for every $0<t<M.$
In particular, if $D\eta(t)\leq 0$ for every $t>0$, we have
$H_{f, \lambda}'(t)>0$ for every $t>0$.
\end{proposition}

\begin{proof}
For $t > 0$, using $f(t)>0$ and $\eta>0$, we get $t f^{\prime}(t)-f(t)>0$.

When $f(t)>\lambda t$, we write
\begin{align}\label{H1}
t H_{f, \lambda}^{\prime}(t) & =\frac{\lambda t}{f(t)}-\left(1-\frac{\lambda t}{f(t)}\right) \frac{t \eta^{\prime}(t)}{\eta(t)^2}\nonumber \\ 
& =\frac{\lambda t}{f(t)}\left[1-\left(\frac{f(t)}{\lambda t}-1\right) \frac{t \eta^{\prime}(t)}{\eta(t)^2}\right] .
\end{align}
If $\chi=0$, then $\eta^{\prime}(t) \leq 0$, and the bracket in (\ref{H1}) is at least one. Hence we have
$H_{f, \lambda}^{\prime}(t)>0$. If $\chi>0$, we further assume $0<t<M$.
Since $f(t) / t$ is strictly increasing, we have
$\frac{f(t)}{\lambda t}<\frac{f(M)}{\lambda M} \leq \hat C .$
Then we get
$$
\begin{aligned}
\left(\frac{f(t)}{\lambda t}-1\right) \frac{t \eta^{\prime}(t)}{\eta(t)^2} & \leq\left(\frac{f(t)}{\lambda t}-1\right) \chi  <( \hat C-1) \chi  \leq 1 .
\end{aligned}
$$
The equation (\ref{H1}) is therefore strictly positive.

When $f(t) \leq \lambda t$, we use the following formula
\begin{align}\label{H2}
t H_{f, \lambda}^{\prime}(t)=\frac{\left[t f^{\prime}(t)-f(t)\right]^2-[f(t)-\lambda t]\left[t^2 f^{\prime \prime}(t)-t f^{\prime}(t)+f(t)\right]}{\left[t f^{\prime}(t)-f(t)\right]^2} .
\end{align}
In this case, the  assumption $D^2\left(\frac{f(t)}{t}\right) \geq 0$ shows
$t^2 f^{\prime \prime}(t)-t f^{\prime}(t)+f(t) \geq 0 .$
Consequently, $H_{f, \lambda}^{\prime}(t)>0 .$
\end{proof}

\begin{lemma}\label{prop:log-structure}
Let $N\ge1$, $0<s<1$, $\lambda>0$, and
\[
f(t)=ct^p\log(1+t^\delta),\quad \text{where} \quad c>0,\quad p\ge1,\quad\delta>0.
\]
Assume $p+\delta<p_s$. Then, $H_{f, \lambda}'(t)>0$ for every $t>0$.
\end{lemma}

\begin{proof}
It suffices to check that $f$ satisfies the conditions in Proposition \ref{thm:general-quotient-monotonicity}, i.e.,
$\eta>0$, $D\eta \le 0$ and $D^2(f/t)>0$ on $(0,\infty)$.

Direct computations yield
$$
f^{\prime}(t)=c p t^{p-1} \log \left(1+t^\delta\right)+\frac{c \delta t^{p+\delta-1}}{1+t^\delta},
$$
and
$$
\eta(t)=p-1+\frac{\delta t^\delta}{\left(1+t^\delta\right) \log \left(1+t^\delta\right)}>0, \quad \forall t > 0.
$$
Moreover,  for $t>0$,
\begin{align*}
t \eta^{\prime}(t) & =\frac{\delta^2 t^\delta}{\left(1+t^\delta\right) \log \left(1+t^\delta\right)}\left[\frac{1}{1+t^\delta}-\frac{t^\delta}{\left(1+t^\delta\right) \log \left(1+t^\delta\right)}\right] \\
& =\frac{\delta^2 t^\delta\left[\log \left(1+t^\delta\right)-t^\delta\right]}{\left(1+t^\delta\right)^2\left[\log \left(1+t^\delta\right)\right]^2}\\
&<0.
\end{align*}
Finally, it is readily seen that
$$D^2\left(\frac{f(t)}{t}\right)=c t^{p-1}[(p-1)^2 \log \left(1+t^\delta\right) +\frac{2(p-1) \delta t^\delta}{1+t^\delta}+\frac{\delta^2 t^\delta}{\left(1+t^\delta\right)^2}]>0.$$
\end{proof}

\section{Uniqueness of the positive ground state solution: proof of Theorems \ref{thm}, \ref{thmN1}}\label{uniqueness}

In this section, we give the proof of Theorems \ref{thm} and Theorems \ref{thmN1} by using a continuation argument. Here, inspired by \cite{SH}, we use a 
homotopic argument by continuation in $\lambda$. The first step is to construct 
a global branch. 

Throughout this section, we always assume that $f$ satisfies the conditions in Proposition \ref{rad}, that is $(f_1)$, $(f_2)$ hold and $H_{f,\lambda}'(t) > 0$ for $0 < t < u(0)$ for any positive, radial, strictly decreasing solution \(u\in H^s(\mathbb R^N)\) of (\ref{Q1})
with Morse index $1$. Additionally, we introduce the following assumptions:
\begin{itemize}
  \item[$(f_3)$]  There exist $q, r\in (2, 2^*_s)$ such that 
  \begin{align}
    0 < q F(t) \leqslant f(t)t   \leqslant r F(t), \quad  \forall t > 0,
  \end{align}
  where $F(t)=\int_0^tf(s)ds$.
  \item[$(f_4)$] For some $a_0>0$ and $r \in (2,2_s^*)$, we have $f(t)= a_0 t^{r-1} +g(t)$,  
  $g(t)\geqslant0$ and there exist $q \in (2,r)$ and $q_0 \in (q, r)$  such that
  \begin{align}
    |r G(t) -g(t)t| \lesssim t^{q}+t^{q_0}, \quad \forall t \geqslant 0,
  \end{align}  
  where $G(t)=\int_0^t g(s) ds$. 
  Without loss of generality, in the following, we always assume $a_0=1$.
\end{itemize}

\begin{proposition} \label{propgb}
  Suppose $(f_3)$. Let $N \geqslant 1, s \in(0,1)$, and $u_0$ be a positive radial solution of \eqref{Q1} with $\lambda = \lambda_0>0$, having Morse index $1$. Then there is a 
  map $u \in C^1([\lambda_0,\infty),H_{\text{rad}}^{s})$ such that $u(\lambda)$ is positive, solves \eqref{Q1} and has Morse index $1$, with $u(\lambda_0) = u_0$.
\end{proposition}

We will prove Proposition \ref{propgb} in  Subsection \ref{subsecgb}.
Thanks to the non-degeneracy result of Theorem \ref{them2}, we can construct a 
local branch stemming from $\lambda = \lambda_0$ by the implicit function theorem. 
Then, in order to extend this branch to $\infty$, we give some results about 
uniform bounds and decay w.r.t. $\lambda$, see details in Subsection \ref{subsecgb}.

Next, we aim to show that the branch constructed by Proposition \ref{propgb} is 
unique. In fact, we define a one-to-one transform $T_\lambda$ from $H^{s}(\mathbb{R}^N)$ 
to $H^{s}(\mathbb{R}^N)$:
\begin{equation}
T_{\lambda}u = \lambda^{-\frac{1}{r-2}}u(\frac{x}{\lambda^{1/(2s)}}).
\end{equation}
Then we get a branch $v \in C^1([\mu_0,\infty),H_{\text{rad}}^{s})$ such that 
$v(\mu):=T_{\mu^{r-2}}u(\mu^{r-2})$ is positive and solves
\begin{equation}\label{v0}
(-\Delta)^s v+v=\mu^{1-r} f(\mu v) \quad \text {in  } \mathbb{R}^N,
\end{equation}
with $v(\mu_0) = v_0 := T_{\lambda_0}u_{\lambda_0}$, where $\mu= 
\lambda^{\frac{1}{r-2}}$. Then, it suffices to prove the uniqueness of the branch 
$v(\mu)$. As $\mu \to \infty$, we will provide the uniform bounds and decay of $v_\mu$, 
and then show that $v_\mu$ converges to the positive, unique (up to 
translations) and non-degenerate ground state solution of the limit equation 
$(-\Delta)^s v_*+v_*- a_\infty\left|v_*\right|^{r-2}v_*=0$ in $\R^N$ for some $ a_\infty >0$. Then by using implicit 
function arguments again, we can prove the following result. 
\begin{proposition} \label{propub}
  Suppose $(f_3)$ and $(f_4)$. Let $N \geqslant 1, s \in(0,1)$,
  and for $i = 1,2$, let $v_i \in C^1([\mu_0,\infty),H_{\text{rad}}^{s})$ 
  be the global branch such that $v_i(\mu) > 0$ solves \eqref{v0} and has Morse index $1$, then 
  $v_1 \equiv v_2$.
\end{proposition}
See Subsection \ref{secunigb} for the proof of Proposition \ref{propub}. Finally in Subsection 
\ref{secpfmr}, we will complete the proof of Theorem \ref{thm} and Theorem \ref{thmN1}.

\subsection{Construction of the global branch} \label{subsecgb}

 Consider:
\begin{align*}
  G: \ & H^s_{\text {rad }}(\mathbb{R}^N) \times [\lambda_0,\infty)\to H^{-s}(\mathbb{R}^N), \\
  &G(u,\lambda)=(-\Delta)^{s} u+ \lambda u -f(u).
\end{align*}
Note that $u$ solves \eqref{Q1} if and only if $G(u,\lambda)=0$. In this subsection, 
we set
$$
L_{+,\lambda} := (-\Delta)^s + \lambda - f'(u).
$$

\begin{lemma}\label{prop4.15} 
Suppose $(f_3)$. Let $u_0 \in H_{\text {rad }}^{s}(\mathbb{R}^N)$ be a positive solution of $G(u_0, \lambda_0)=0$ having Morse index $1$.   
Then for some $\delta>0$, there exists a map 
$u(\lambda) \in C^1(\left[\lambda_0,\lambda_0+\delta),H_{\text{rad}}^{s}\right)$ 
such that the following holds true,  where we denote $u_\lambda=u(\lambda)$ in the sequel.

(i) $u_\lambda$ solves $G(u, \lambda)=0$ and has Morse index $1$ for all $\lambda \in [\lambda_0,\lambda_0+\delta)$.

(ii) There exists $\epsilon>0$ such that $u_\lambda$ is the unique solution of 
$G(u, \lambda)=0$ for $\lambda \in[\lambda_0,\lambda_0+\delta)$ in the 
neighborhood $\left\{u \in H_{\text {rad }}^{s}:\left\|u-u_\lambda
\right\|_{H^{s}}<\epsilon\right\}$. In particular, we have that $u_{\lambda_0}=u_0$ holds.
\end{lemma}
\begin{proof} 
We use an implicit function argument for the map
$$
\Psi: H_{\text {rad }}^{s} \times [\lambda_0,\lambda_0+\delta) \rightarrow 
H_{\text {rad }}^{s}, \quad \Psi(u, \lambda)=u-\frac{1}{(-\Delta)^s+\lambda} f( u).
$$
It is not difficult to see that $\Psi$ is a well-defined map of class $C^1$. 
Furthermore, $\Psi\left(u_0, \lambda_0\right)=0$ and $\Psi(u, \lambda)=0$ if 
and only if $u \in H_{\text {rad }}^{s}$ solves $G(u,\lambda)=0$. Next, we consider the Fréchet derivative
$$
\partial _u \Psi\left(u_0, \lambda_0\right)=1-K, K=\frac{1}{(-\Delta)^s+\lambda_0} f'(u_0)
$$
By  Theorem \ref{them2}, $0$ is not an eigenvalue of $L_{+,\lambda_0}=
(-\Delta)^s+\lambda_0-f'(u_0)$ constrained on $L_{rad }^2$. Note that the 
operator $K$ is compact on $L_{\text {rad }}^2$ and we have that 
$1 \notin \sigma(K)$. Since $f'(u_0) \varphi \in H^{-s}$ if $\varphi \in 
H_{\text {rad }}^{s}$, we have that $K$ maps $H_{\text {rad }}^{s}$ to 
$H_{\text {rad }}^{s}$. Thus the bounded inverse $(1-K)^{-1}$ exists on 
$H_{\text {rad }}^{s}$, i.e., the Fréchet derivative $\partial_u\Psi$ has a 
bounded inverse on $H_{\text {rad }}^{s}$ at $\left(u_0, \lambda_0 \right)$. 
By the implicit function theorem, we deduce that the assertions 
$(i)$ and $(i i)$ in Lemma \ref{prop4.15} hold for some $\delta>0$ and 
$\epsilon>0$ sufficiently small.
\end{proof}
We now turn to the global extension of the locally unique branch 
$u_\lambda$ constructed in the above Lemma \ref{prop4.15}. First, we will give 
some uniform a priori estimates. We write $a \lesssim b$ to denote that $a\leqslant C b$ with some constant $C>0$ that may change from line to line and only depends on $N$, $ s, \lambda$ and $u_0$. Furthermore, we use $a \sim b$ to denote that both $a \lesssim b$ and $b \lesssim a$ hold.
In what follows,  let us assume that 
$N > 2s$ holds. The arguments can be adapted with some modifications to the case $N \leqslant 2s$.

\begin{lemma}[Pohozaev identity, \cite{SH}, Lemma 5.6]\label{Pohozaev}   
Let $u$ be a weak solution of \eqref{Q1},  then the following integral identity holds true:

\begin{equation}\label{a4.4}
(N-2s)\int_{\mathbb{R}^N}  |(-\Delta)^{s/2} u|^2dx+ N\lambda \int_{\mathbb{R}^N} u^2 dx= 
2N\int_{\mathbb{R}^N}F(u) dx.
\end{equation}
\end{lemma}

\begin{lemma}\label{lem4.16}
Suppose $(f_3)$. Then either $\lambda_* = \infty$ or we have the following power estimates
\[
\int_{\mathbb{R}^N} f(u_{\lambda})u_{\lambda} dx  \sim \lambda 
\int_{\mathbb{R}^N}  u_{\lambda}^2dx   \sim  \int_{\mathbb{R}^N}  |(-\Delta)^{s/2} u_{\lambda}|^2dx \sim 1,
\]
for all $\lambda\in [\lambda_0, \lambda_*)$, where $u_{\lambda}$ is given in Lemma \ref{prop4.15} and
\begin{align*}
  \lambda_*= \sup \{\hat \lambda >0: \ & u_{\lambda}\in C^1([\lambda_0,\hat\lambda),
H^s_{\text {rad }}(\mathbb{R}^N)), \\
& u_{\lambda}\text{ satisfies the assumptions of Lemma \ref{prop4.15}}\text{ for} \lambda\in [\lambda_0, \hat \lambda)\}.
\end{align*}
\end{lemma}
\begin{proof}
	We assume $\lambda_* < \infty$ and prove the estimates stated in the lemma.

\vskip0.1in
{\bf Step 1.}  We claim $\int_{\mathbb{R}^N} f(u_{\lambda})u_{\lambda} dx\gtrsim 1$ for $\lambda\in (\lambda_0,\lambda_*)$.

By integrating $ G(u_{\lambda},\lambda)=0$ after a multiplication by 
$u_{\lambda}$, we have
\begin{equation}\label{a4.21}
\int_{\mathbb{R}^N}  |(-\Delta)^{s/2} u_{\lambda}|^2dx + \lambda \int_{\mathbb{R}^N} u_{\lambda}^2 dx= \int_{\mathbb{R}^N} f(u_{\lambda})u_{\lambda} dx.
\end{equation}
Consequently, (\ref{a4.4}) and (\ref{a4.21}) together, give

\begin{equation}\label{a4.22}
2s \lambda \int_{\mathbb{R}^N}  u_{\lambda}^2dx= 2\int_{\mathbb{R}^N}  N F(u_{\lambda})dx -(N-2s) \int_{\mathbb{R}^N} f(u_{\lambda})u_{\lambda} dx,
\end{equation}
\begin{equation}\label{a4.23}
2s \int_{\mathbb{R}^N}|(-\Delta)^{s/2} u_{\lambda}|^2dx = 
N \int_{\mathbb{R}^N} f(u_{\lambda})u_{\lambda} dx - 
2N\int_{\mathbb{R}^N}F(u_{\lambda})dx.
\end{equation}
By  assumption $(f_3)$,   we have  $ F(t)  \sim f(t)t$, therefore we derive 
\begin{equation}\label{MV}
\int_{\mathbb{R}^N} f(u_{\lambda})u_{\lambda} dx  \sim \lambda \int_{\mathbb{R}^N}  u_{\lambda}^2dx   \sim  \int_{\mathbb{R}^N}  |(-\Delta)^{s/2} u_{\lambda}|^2dx.
\end{equation}
Since $2<q<r<2^{*}_{s}$, $\forall u\in H^s(\mathbb{R}^N),$
from the following fractional Gagliardo-Nirenberg-Sobolev inequality
\begin{equation}\label{GNS}
\int_{\mathbb{R}^N}|u|^{q}dx\leqslant C_{opt}\left(\int_{\mathbb{R}^N}|(-\Delta)^{s/2}u|^2dx\right)^{\frac{N(q-2)}{4s}}
\left(\int_{\mathbb{R}^N}|u|^2dx\right)^{\frac q2-\frac{N(q-2)}{4s}}, 
\end{equation}
where $C_{opt}$ is the best constant and depends on $s, N, q$, 
we   derive

\begin{eqnarray}\label{a4.27}
\int_{\mathbb{R}^N} f(u_{\lambda})u_{\lambda} dx &\lesssim & \int_{\mathbb{R}^N} u_{\lambda}^q+ u_{\lambda}^r dx \nonumber \\
& \lesssim & \left(\int_{\mathbb{R}^N}|(-\Delta)^{s/2}u|^2dx\right)^{\frac{N(q-2)}{4s}}
\left(\int_{\mathbb{R}^N}|u|^2dx\right)^{\frac q2-\frac{N(q-2)}{4s}}  \nonumber \\
&& +\left (\int_{\mathbb{R}^N}|(-\Delta)^{s/2}u|^2dx\right)^{\frac{N(r-2)}{4s}}
\left(\int_{\mathbb{R}^N}|u|^2dx\right)^{\frac r2-\frac{N(r-2)}{4s}}
 \nonumber \\
&\lesssim &  \left( \int_{\mathbb{R}^N} f(u_{\lambda})u_{\lambda} dx \right)^{q/2} \lambda^{  \frac{N(q-2)}{4s}-\frac q2 } \nonumber \\
&& + \left( \int_{\mathbb{R}^N} f(u_{\lambda})u_{\lambda} dx \right)^{r/2} \lambda^{  \frac{N(r-2)}{4s}-\frac r2 } .
\end{eqnarray}
Thus $\int_{\mathbb{R}^N} f(u_{\lambda})u_{\lambda} dx\gtrsim 1$ for $\lambda\in (\lambda_0,\lambda_*)$.

\vskip0.1in
{\bf Step 2.}
Next we show that $\int_{\mathbb{R}^N} f(u_{\lambda})u_{\lambda} dx\lesssim 1$ for $\lambda\in (\lambda_0,\lambda_*)$.\\
Set
\[g(\lambda)=\int_{\mathbb{R}^N}[f(u_{\lambda})u_{\lambda} -2 F(u_{\lambda})]dx.
\]
Assumption $(f_3)$ and (\ref{MV}), yield  that 
\begin{equation}\label{fV}
\lambda \int_{\mathbb{R}^N} u^2_{\lambda} dx \sim 
\int_{\mathbb{R}^N} f(u_{\lambda})u_{\lambda} dx  \sim g (\lambda).
\end{equation}
Therefore, it suffices to show that  $g (\lambda)$ has uniform upper bound  on $ (\lambda_0,\lambda_*)$.\\
Differentiating $G(u_{\lambda},\lambda)=0$ with respect to $\lambda$, we get
\begin{equation}\label{02.8}
(-\Delta)^s \frac{d}{d\lambda} u_{\lambda}=-  u_{\lambda} -\lambda \frac{d}{d\lambda} u_{\lambda}+ f'( u_{\lambda})\frac{d}{d\lambda} u_{\lambda},
\end{equation}
hence 
\begin{equation}\label{02.800}
L_{+,\lambda}(\frac{d}{d\lambda} u_{\lambda})=-  u_{\lambda}.
\end{equation}
Note that
$L_{+,\lambda} u_{\lambda}=f(u_{\lambda})- f'(u_{\lambda})u_{\lambda},$
and this implies, via  (\ref{02.800}),
\begin{eqnarray}\label{gg}
g' (\lambda) &= & \int_{\mathbb{R}^N} f'( u_{\lambda})\frac{d}{d\lambda} u_{\lambda} \cdot  u_{\lambda}  + f( u_{\lambda})\frac{d}{d\lambda} u_{\lambda} - 2f(u_{\lambda}) \frac{d}{d\lambda}u_{\lambda}  dx\nonumber \\
&= & \int_{\mathbb{R}^N}  f'( u_{\lambda})u_{\lambda} \frac{d}{d\lambda} u_{\lambda}-f( u_{\lambda}) \frac{d}{d\lambda} u_{\lambda} dx \nonumber \\
&= & -\langle L_{+,\lambda} u_{\lambda}, \frac{d}{d\lambda} u_{\lambda}\rangle _{L^2} \nonumber \\
&= & -\langle  u_{\lambda}, L_{+,\lambda} \frac{d}{d\lambda} u_{\lambda}\rangle _{L^2} \nonumber \\
&=&\int_{\mathbb{R}^N} u_{\lambda}^2 dx.
\end{eqnarray}
Via (\ref{fV}) and (\ref{gg}), we obtain the differential inequality 
\begin{eqnarray}\label{818}
g' (\lambda)&  \lesssim &  \frac{1}{\lambda }g(\lambda),
\end{eqnarray}
which by integration yields the uniform upper bound  $g(\lambda) \lesssim 1$  for  $\lambda\in [\lambda_0, \lambda_*)$. 
The proof of Lemma \ref{lem4.16}  is now complete. 
\end{proof}

\begin{lemma}\label{lem4.17}
Suppose $(f_3)$ and that $u_{\lambda_0}(x)>0$ is  ground state. Let $\lambda_*$ is 
given by Lemma \ref{lem4.16}.
Then, for all $\lambda \in\left[\lambda_0,\lambda_*\right)$, we have
\begin{itemize}
  \item[(i)] $u_{\lambda}(x)>0$ for all $x\in \mathbb{R}^N$;
  \item[(ii)] (uniform decay estimate) either $\lambda_* = \infty$ or $u_{\lambda}(x)\leqslant C |x|^{-(N+2s)}$ for 
  $|x|> R$, where $R$ is some constant independent of $\lambda$.
\end{itemize}
\end{lemma}
\begin{proof}
Define the linearized operator

\begin{equation}\label{L-}
L_{-,\lambda}=(-\Delta)^s +\lambda - \frac{f(u)}{u}  
\end{equation}
acting on $L^2(\mathbb{R}^N)$.
Obviously, if $\lambda\to \tilde \lambda \in [\lambda_0,\lambda_*)$, then $ L_{\pm,\lambda}\to L_{\pm,\tilde \lambda}$ in norm-resolvent sense.

The proof of $(i)$ is divided into two steps. First, we show that if $u_{\tilde \lambda}>0$ for some
$\tilde \lambda \in [\lambda_0,\lambda_*)$, then $u_{\lambda}>0$ for $\lambda$ 
close to $\tilde \lambda$. Note that $L_{-,\lambda} u_{\lambda}=0$, i.e., $0$ 
is an eigenvalue of $L_{-,\lambda} $ and  $u_{\lambda}$ is an eigenfunction with 
respect to $0$.  Since $\sigma_{ess}( L_{-,\lambda})= [\lambda, +\infty)$, we 
know that  $\inf\sigma( L_{-,\lambda})$ is an eigenvalue.  
\cite[Lemma C.4]{FS} yields that the first eigenvalue is simple and the 
corresponding eigenfunction can be chosen strictly positive. Thus
$\lambda_1(L_{-,\tilde \lambda})=0$ is simple and we choose $\psi_{1,\tilde \lambda}=u_{\tilde \lambda}>0$ as the corresponding eigenfunction, where $\lambda_1(L_{-,\tilde \lambda})$  means the first eigenvalue of  $L_{-,\tilde \lambda}$. Since  $L_{-, \lambda}\to L_{-,\tilde \lambda}$
in norm-resolvent sense when $\lambda\to \tilde \lambda\in 
[\lambda_0,\lambda_*)$, we derive that $ \lambda_1(L_{-,\lambda})\to 0$. 
Since $\lambda_1(L_{-,\tilde \lambda})$ is simple and isolated, we assume  
that $\sigma (L_{-, \lambda})\cap (-\varepsilon, \varepsilon)=
\{\lambda_1(L_{-,\lambda})\} $ for some small $\varepsilon$ and  
$\lambda$ close enough to $\tilde \lambda$. Furthermore, $0$ is an eigenvalue of 
$L_{-,\lambda}$, showing $\lambda_1(L_{-, \lambda})=0$. Since $ u_{\lambda} 
\to u_{\tilde \lambda}>0$ in $H^{s}$, we obtain $u_{ \lambda}>0$ for $\lambda$ close enough to $\tilde \lambda$.

The second step is to show that if $\{ \lambda_n \}_{n=1}^{+\infty} \subset [\lambda_0,\lambda_*)$  is a sequence with  $\lambda_n\to \hat \lambda$ such that $ u_{ \lambda_n}>0$ for some $\hat \lambda\in
[\lambda_0,\lambda_*)$, then  $ u_{\hat \lambda}>0$. Obviously, $ u_{\hat \lambda}
\geqslant 0$ since $ u_{ \lambda_n} \to u_{\hat \lambda}$ in  $H^{s}$. Thus $0$ 
is  the first eigenvalue of  $  L_{-, \hat \lambda}$ and so $u_{\hat \lambda} > 0$. 
The proof of $(i)$ is complete.

Next, we prove $(ii)$. Since $u_{ \lambda}$ is monotonically decreasing with 
respect to $r=|x|$, we have $|x|^N  u^2_{ \lambda}(x)
\lesssim u^2_{ \lambda}(x)\int_{B_{|x|}} dy \lesssim 
\int_{B_{|x|}} u^2_{\lambda}(y) dy \lesssim 1.$ We use the uniform a priori 
estimate shown in Lemma \ref{lem4.16} in the final inequality. Thus $u_{ \lambda}(x)\lesssim  |x|^{-\frac N2}$.

We write $u_{ \lambda}=\frac{1}{(-\Delta)^s+ \theta}z_{\lambda}$ where 
$0< \theta< \lambda$ and $z_{\lambda}=(\theta-\lambda)u_{\lambda}+
f(u_{ \lambda})$. Since $u_{ \lambda}(x)\lesssim  |x|^{-\frac N2}$, 
we can take $R$ large enough such that $supp z_{\lambda}^+\subset B_{\frac R2 }$ 
and $R$ is independent of $\lambda$.  Suppose $G_{s,\theta}(x)$ is the kernel of 
the resolvent $(( -\Delta)^s+\theta )^{-1}$ on $ \mathbb{R}^N$.
Notice that  $f(u)\lesssim (u^{q-1}+u^{r-1})$, then we obtain
\begin{eqnarray}\label{a4.32}
\int_{B_{\frac R2}} z_{\lambda}^+(y)dy
&\lesssim&\int_{B_{\frac R2} \cap supp z_{\lambda}^+  }  (u_{\lambda}^{q-1}+u^{r-1}_{\lambda}) dy \nonumber \\
&\lesssim & (\int_{B_{\frac R2}  }1 dy  )^{\frac{1}q}
(\int_{B_{\frac R2} \cap supp z_{\lambda}^+  } u_{\lambda}^{q} dy  )^{\frac{q-1}q}+ (\int_{B_{\frac R2}  }1 dy  )^{\frac{1}r}
(\int_{B_{\frac R2} \cap supp z_{\lambda}^+  }u^{r}_{\lambda} dy)^{\frac{r-1}r} \nonumber \\
&\lesssim &  R^{\frac Nq} +R^{\frac Nr}
\end{eqnarray}
We use the uniform a priori estimates shown in Lemma \ref{lem4.16} in the final inequality. 
Then, we can conclude that: for all $|x|>R$,
\begin{eqnarray}\label{a4.33}
u_{\lambda}(x)
&=&\int_{\mathbb{R}^N} G_{s,\theta}(x-y)z_{\lambda}(y)dy \leqslant\int_{B_{\frac R2}} G_{s,\theta}(x-y)z^+_{\lambda}(y)dy \nonumber \\
&\lesssim & \int_{B_{\frac R2}}|x-y|^{-(N+2s)} z^+_{\lambda}(y)dy \lesssim \int_{B_{\frac R2}}|x|^{-(N+2s)} z^+_{\lambda}(y)dy \nonumber \\
&\lesssim & |x|^{-(N+2s)}.
\end{eqnarray}
Hence, we complete the proof.
\end{proof}

\begin{lemma}\label{lem4.18}
Suppose $(f_3)$. If $\lambda_* < \infty$,
let $\{\lambda_n\}_{n=1}^{+\infty}\subset 
[\lambda_0,\lambda_*)$ be a sequence such that $\lambda_n\to \lambda_*$. Then, after possibly passing to a 
subsequence, we have $u_{\lambda_n}\to u_{\lambda_*}$ in 
$H^s_{\text {rad }}(\mathbb{R}^N)$, where $u_{\lambda_*}(|x|)>0$ satisfies
\begin{equation}\label{a4.34}
(-\Delta)^s u_{\lambda_*}+\lambda_* u_{\lambda_*}= f(u_{\lambda_*})  \quad \text{in } \mathbb{R}^N,
\end{equation}
and $\mu(u_{\lambda_*}) = 1$.
\end{lemma}

\begin{proof}
By Lemma \ref{lem4.16}, $u_{\lambda_n}$ is bounded in $H^s$. Passing to a 
subsequence if necessary, we assume that $u_{\lambda_n} \rightharpoonup  u_{\lambda_*}$ 
in $H^s, u_{\lambda_n} \rightarrow u_{\lambda_*}$ in $L_{loc}^p$ for 
$1 \leqslant p<2_s^*, u_{\lambda_n} \rightarrow u_{\lambda_*}$ a.e. on $\mathbb{R}^N$. 
With the help of the uniform decay estimate shown in Lemma \ref{lem4.17} $(ii)$, 
we obtain that $u_{\lambda_n} \rightarrow u_{\lambda_*}$ in $L^p$ for 
$1\leqslant p<2_s^*$. Bootstrap arguments implies that $u_{\lambda_n} \rightarrow 
u_{\lambda_*}$ in $H^{s}$. Obviously, $u_{\lambda_*}\geqslant 0$ is nonzero and satisfies equation 
(\ref{a4.34}). Thus $0$ is the first eigenvalue of $L_{-, \lambda_*}$, showing $u_{\lambda_*}>0$. By the lower semi-continuity of the Morse index, we have $\mu(u_{\lambda_*}) \le 1$. On the other hand, using $\langle u_{\lambda_*},L_{+, \lambda_*} u_{\lambda_*}\rangle_{L^2} <0$ we deduce that $\mu(u_{\lambda_*}) \ge 1$. Thus $\mu(u_{\lambda_*}) = 1$.
\end{proof}

\begin{proof}[Proof of Proposition \ref{propgb}]
  It suffices to prove that $\lambda_* = \infty$, which is given by Lemma \ref{lem4.16}. 
  Suppose on the contrary that $\lambda_* < \infty$. By Lemma \ref{lem4.17}, we can assume that 
  $u(\lambda)$ is positive and has Morse index $1$ for all $\lambda \in [\lambda_0,\lambda_*)$. Then, by Lemma \ref{lem4.18}, 
  we get $u_{\lambda_*}$ which is a positive solution of \eqref{a4.34} with $\mu(u_{\lambda_*}) = 1$. Applying Theorem \ref{them2}, 
  we obtain that $u_{\lambda_*}$ is non-degenerate. Then using implicit function arguments again 
  (like the proof of Lemma \ref{prop4.15}), we can extend the branch $u(\lambda)$ to $\lambda_* + \delta$, 
  contradicting the definition of $\lambda_*$. Hence, we complete the proof.
\end{proof}

\subsection{Uniqueness of the global branch} \label{secunigb}

In this subsection, we aim to prove Proposition \ref{propub}, and, 
we set
$$
L_{+,\mu}=(-\Delta)^s +1 -\mu^{2-r} f'(\mu v).
$$

\begin{lemma}\label{lem8.1}
Suppose $(f_3)$ and $(f_4)$. Let $v(\mu) \in C^1([\mu_0,\infty),H_{\text{rad}}^{s})$ 
be the global branch such that $v_\mu := v(\mu)$ is positive and solves \eqref{v0}. 
Then we have the a priori bounds
\[
\mu^{1-r}\int_{\mathbb{R}^N} f(\mu v_{\mu})v_{\mu} dx   \sim  \int_{\mathbb{R}^N} v_{\mu}^2dx   \sim  \int_{\mathbb{R}^N}  |(-\Delta)^{s/2}v_{\mu}|^2dx \sim 1,
\]
uniformly w.r.t. $\mu\in [\mu_0,\infty)$.
\end{lemma}
\begin{proof}
 We divide the proof into the following two steps.

 \vskip0.1in
{\bf Step 1.} Lower bounds.  

First, we claim 
\[
\mu^{1-r}\int_{\mathbb{R}^N} f(\mu v_{\mu})v_{\mu} dx  \sim \int_{\mathbb{R}^N}  v_{\mu}^2dx   \sim  \int_{\mathbb{R}^N}  |(-\Delta)^{s/2} v_{\mu}|^2dx.
\]
By integrating \eqref{v0} after a multiplication by $v_{\mu}$, we obtain
\begin{equation}\label{v4.21}
\int_{\mathbb{R}^N}  |(-\Delta)^{s/2} v_{\mu}|^2dx +  \int_{\mathbb{R}^N} v_{\mu}^2 dx= \int_{\mathbb{R}^N}\mu^{1-r} f(\mu v_{\mu})v_{\mu} dx.
\end{equation}
For \eqref{v0}, we have the  Pohozaev type identity
\begin{equation}\label{ph}
(N-2s)\int_{\mathbb{R}^N}  |(-\Delta)^{s/2}v_{\mu}|^2dx+N \int_{\mathbb{R}^N} v_{\mu}^2 dx= 
2N\int_{\mathbb{R}^N} \mu^{-r} F(\mu v_{\mu}) dx.
\end{equation}
Via  (\ref{v4.21})   and  (\ref{ph}), we have
\begin{equation}\label{v4.22}
2s \int_{\mathbb{R}^N} v_{\mu} ^2dx= 2N\int_{\mathbb{R}^N}\mu^{-r} F(\mu v_{\mu})dx -(N-2s) \int_{\mathbb{R}^N}\mu^{1-r} f(\mu v_{\mu})v_{\mu} dx,
\end{equation}
and
\begin{equation}\label{v4.23}
2s \int_{\mathbb{R}^N}|(-\Delta)^{s/2} v_{\mu} |^2dx = 
N\mu^{1-r}\int_{\mathbb{R}^N}f(\mu v_{\mu})v_{\mu}dx - 2N\mu^{-r} \int_{\mathbb{R}^N}F(\mu v_{\mu} ) dx.
\end{equation}
Assumption $(f_3)$ implies  $F(t)\sim f(t)t$, hence 
\[
\mu^{1-r}\int_{\mathbb{R}^N} f(\mu v_{\mu})v_{\mu} dx  \sim \int_{\mathbb{R}^N}  v_{\mu}^2dx   \sim  \int_{\mathbb{R}^N}  |(-\Delta)^{s/2} v_{\mu}|^2dx.\]
Since $2<q<r<2^*_s$,
from the  fractional Gagliardo-Nirenberg-Sobolev inequality (\ref{GNS}),
 we  can derive

\begin{eqnarray}\label{v4.27}
\mu^{1-r} \int_{\mathbb{R}^N} f(\mu v_{\mu})v_{\mu} dx &\lesssim &    \int_{\mathbb{R}^N}\mu^{q-r} v_{\mu}^q +\mu^{r-r} v_{\mu}^r dx \nonumber \\
& \lesssim & \mu^{q-r} (\int_{\mathbb{R}^N}|(-\Delta)^{s/2}v_{\mu}|^2dx)^{\frac{N(q-2)}{4s}}
\left(\int_{\mathbb{R}^N}|v_{\mu}|^2dx\right)^{\frac q2-\frac{N(q-2)}{4s}}  \nonumber \\
&& + (\int_{\mathbb{R}^N}|(-\Delta)^{s/2}v_{\mu}|^2dx)^{\frac{N(r-2)}{4s}}
\left(\int_{\mathbb{R}^N}|v_{\mu}|^2dx\right)^{\frac r2-\frac{N(r-2)}{4s}}
 \nonumber \\
&\lesssim & \mu^{q-r} ( \mu^{1-r} \int_{\mathbb{R}^N} f(\mu v_{\mu})v_{\mu} dx )^{q/2}  \nonumber \\
&& + ( \mu^{1-r} \int_{\mathbb{R}^N} f(\mu v_{\mu})v_{\mu} dx )^{r/2}.
\end{eqnarray}
Thus $\mu^{1-r} \int_{\mathbb{R}^N} f(\mu v_{\mu})v_{\mu} dx\gtrsim 1$ for 
$\mu\in [\mu_0, \infty)$.

\vskip0.1in
\textbf{Step 2.} Upper bounds.  Next we show that $\mu^{1-r} \int_{\mathbb{R}^N} f(\mu v_{\mu})v_{\mu} dx\lesssim 1$ for 
$\mu\in [\mu_0,  \infty)$. 

Set
\[k(\mu)=\mu^{1-r} \int_{\mathbb{R}^N}  f(\mu v_{\mu})v_{\mu}dx,
\]

\[h(\mu)=2\mu^{-r} \int_{\mathbb{R}^N}  F(\mu v_{\mu})dx,
\]
and
\[B(\mu)=\mu^{-r} \int_{\mathbb{R}^N}  f(\mu v_{\mu})\mu v_{\mu}-  2 F(\mu v_{\mu})dx.
\]
By $(f_3)$, we have $B(\mu)\sim k(\mu)$. Now it suffices to prove that $B(\mu) \lesssim 1$.

Clearly, a brief computation gives 
\begin{eqnarray}\label{g2.9}
k' (\mu) &= & (1-r) \mu^{-r}   \int_{\mathbb{R}^N} f(\mu v_{\mu})v_{\mu}dx
+  \mu^{1-r}   \int_{\mathbb{R}^N} f'(\mu v_{\mu}) [ v_{\mu} + \mu \frac{d}{d\mu}  v_{\mu} ] v_{\mu} dx \nonumber \\
&+&  \mu^{1-r}   \int_{\mathbb{R}^N} f(\mu v_{\mu}) \frac{d}{d\mu}  v_{\mu} dx \nonumber \\
&= & (1-r) \mu^{-r}   \int_{\mathbb{R}^N} f(\mu v_{\mu})v_{\mu}dx
+  \mu^{1-r}   \int_{\mathbb{R}^N} f'(\mu v_{\mu})  v_{\mu}^2
+f'(\mu v_{\mu}) \mu  v_{\mu}\frac{d}{d\mu}  v_{\mu}  dx \nonumber \\
&+&  \mu^{1-r}   \int_{\mathbb{R}^N} f(\mu v_{\mu}) \frac{d}{d\mu}  v_{\mu}dx
\end{eqnarray}
and
\begin{eqnarray}\label{h2.9}
h' (\mu) &= & 2(-r) \mu^{-(r+1)}   \int_{\mathbb{R}^N}F(\mu v_{\mu})dx + 2\mu^{-r}\int_{\mathbb{R}^N}
f(\mu v_{\mu}) [ v_{\mu} + \mu \frac{d}{d\mu}  v_{\mu} ]dx.
\end{eqnarray}
Differentiating \eqref{v0} with respect to $\mu$, we have
\begin{equation}\label{v2.8}
(-\Delta)^s \frac{d}{d\mu} v_{\mu}+ \frac{d}{d\mu} v_{\mu}= (1-r) \mu^{-r}  f(\mu v_{\mu})
+ \mu^{1-r}  f'(\mu v_{\mu}) [  v_{\mu}+ \mu \frac{d}{d\mu} v_{\mu} ],
\end{equation}
consequently,
\begin{equation}\label{dv}
L_{+,\mu}(\frac{d}{d\mu} v_{\mu} )= \mu^{-r} \left[   f'(\mu v_{\mu})  \mu v_{\mu} + (1-r)   f(\mu v_{\mu}) \right].
\end{equation}
Note that
$L_{+,\mu}  v_{\mu}=\mu^{1-r} [f(\mu v_{\mu}) -  f'(\mu v_{\mu})\mu v_{\mu} ],$  then we go further
\begin{eqnarray}\label{B2.9}
k' (\mu)- h' (\mu)
&= & (1-r) \mu^{-r}   \int_{\mathbb{R}^N} f(\mu v_{\mu})v_{\mu}dx
+  \mu^{1-r}   \int_{\mathbb{R}^N} f'(\mu v_{\mu})  v_{\mu}^2
+f'(\mu v_{\mu}) \mu  v_{\mu}\frac{d}{d\mu}  v_{\mu}  dx \nonumber \\
&& +  \mu^{1-r}   \int_{\mathbb{R}^N} f(\mu v_{\mu}) \frac{d}{d\mu}  v_{\mu} dx\nonumber \\
&& - 2(-r) \mu^{-(r+1)}   \int_{\mathbb{R}^N}F(\mu v_{\mu})dx - 2\mu^{-r}\int_{\mathbb{R}^N}
f(\mu v_{\mu}) [ v_{\mu} + \mu \frac{d}{d\mu}  v_{\mu} ]  dx\nonumber \\
&= & (1-r) \mu^{-r}   \int_{\mathbb{R}^N} f(\mu v_{\mu})v_{\mu}dx
+  \mu^{1-r}   \int_{\mathbb{R}^N} f'(\mu v_{\mu})  v_{\mu}^2  dx \nonumber \\
&& + \mu^{1-r}   \int_{\mathbb{R}^N} f'(\mu v_{\mu}) \mu  v_{\mu}\frac{d}{d\mu}  v_{\mu}  dx
-  \mu^{1-r}   \int_{\mathbb{R}^N} f(\mu v_{\mu}) \frac{d}{d\mu}  v_{\mu} dx\nonumber \\
&& - 2(-r) \mu^{-(r+1)}   \int_{\mathbb{R}^N}F(\mu v_{\mu})dx - 2\mu^{-r}\int_{\mathbb{R}^N}
f(\mu v_{\mu})  v_{\mu} dx  \nonumber \\
&= &  \mu^{-(r+1)}   \int_{\mathbb{R}^N} (1-r) f(\mu v_{\mu})\mu v_{\mu}
+  f'(\mu v_{\mu}) (\mu v_{\mu})^2  dx \nonumber \\
&& - \langle L_{+,\mu}  v_{\mu}, \frac{d}{d\mu}  v_{\mu} \rangle_{L^2}
+ 2 \mu^{-(r+1)}   \int_{\mathbb{R}^N}r F(\mu v_{\mu}) -
f(\mu v_{\mu}) \mu v_{\mu} dx. 
\end{eqnarray}
Clearly,
\begin{eqnarray}\label{Lvv}
\langle L_{+,\mu}  v_{\mu}, \frac{d}{d\mu}  v_{\mu} \rangle_{L^2} & =& \langle   v_{\mu},L_{+,\mu} \frac{d}{d\mu}  v_{\mu} \rangle _{L^2} \nonumber \\
& =&  \int_{\mathbb{R}^N}  \mu^{-r} v_{\mu} \left[   f'(\mu v_{\mu})  \mu v_{\mu} + (1-r)   f(\mu v_{\mu}) \right] dx,
\end{eqnarray}
recalling that  (\ref{dv}) holds. Going back to (\ref{B2.9}), we have
\begin{eqnarray}\label{vg}
B' (\mu)
&= &  \mu^{-(r+1)}   \int_{\mathbb{R}^N} (1-r) f(\mu v_{\mu})\mu v_{\mu}
+  f'(\mu v_{\mu}) (\mu v_{\mu})^2  dx \nonumber \\
&-& \langle L_{+,\mu}  v_{\mu}, \frac{d}{d\mu}  v_{\mu} \rangle_{L^2}
 + 2 \mu^{-(r+1)}   \int_{\mathbb{R}^N}r F(\mu v_{\mu}) -
f(\mu v_{\mu}) \mu v_{\mu} dx  \nonumber \\
&= & 2 \mu^{-(r+1)}   \int_{\mathbb{R}^N}r F(\mu v_{\mu}) -
f(\mu v_{\mu}) \mu v_{\mu} dx.
\end{eqnarray}
It is easy to see $B(\mu)>0, B'(\mu)\ge0$.

Now, we are ready to come to the main point of the proof of Lemma \ref{lem8.1}. 
It is convenient to use the following notation
$$
M_{\mu}=\int_{\mathbb{R}^N}\left| v_{\mu} \right|^2, \quad T_{\mu}=\int_{\mathbb{R}^N}\left|(-\Delta)^{s / 2} v_{\mu} \right|^2, \quad V_{\mu}=\mu^{-r} \int_{\mathbb{R}^N}  f(\mu v_{\mu})\mu v_{\mu},
$$
for  $\mu\in [\mu_0, \infty)$. We first recall that
$$
M_{\mu} \sim T_{\mu} \sim V_{\mu}.
$$

Note that
\begin{equation}\label{fmu}
\mu^{-r} \int_{\mathbb{R}^N}  f(\mu v_{\mu})\mu v_{\mu}=
\int_{\mathbb{R}^N}v_{\mu}^r dx+\mu^{-r}\int_{\mathbb{R}^N}g(\mu v_{\mu})\mu v_{\mu}dx, 
\end{equation}
by $(f_4)$. For any $p\in (2,r)$,   in light of  (\ref{fmu}),  we have
\begin{eqnarray}\label{ine}
\int_{\mathbb{R}^N}|v_{\mu}|^p dx
&\lesssim & (\int_{\mathbb{R}^N}|v_{\mu}|^2 dx)^{\frac {\alpha p}2} (\int_{\mathbb{R}^N}|v_{\mu}|^r dx)^{\frac {(1-\alpha) p}r}\nonumber \\
&\lesssim & (M_\mu)^{\frac {\alpha p}2} (V_\mu)^{\frac {(1-\alpha) p}r}\nonumber \\
&\sim & M_\mu,
\end{eqnarray}
where the last inequality holds by  $\frac 1p =\frac {\alpha}2+\frac {1-\alpha}r$. 
Then by $(f_4)$ again, we derive that
\begin{eqnarray}\label{B'1}
B'(\mu) &= & 2 \mu^{-(r+1)}   \int_{\mathbb{R}^N}r F(\mu v_{\mu}) -
f(\mu v_{\mu}) \mu v_{\mu} dx \nonumber \\
&= & 2 \mu^{-(r+1)}  \int_{\mathbb{R}^N}  rG(\mu v_{\mu}) -g( \mu v_{\mu})\mu v_{\mu} dx\nonumber \\
&\lesssim&  \mu^{-(r+1)}  \int_{\mathbb{R}^N} (\mu v_{\mu})^q +( \mu v_{\mu})^{q_0} dx   \nonumber \\
&\lesssim &   \mu^{-(r-q+1)}   M_{\mu}  +  \mu^{-(r-q_0+1)}M_{\mu},
\end{eqnarray}
where the last inequality follows from (\ref{ine}). Hence
\begin{eqnarray}\label{Bb'}
\frac{B'(\mu)}{B(\mu)}
&\lesssim & \frac{   \mu^{-(r-q+1)}   M_{\mu}  +  \mu^{-(r-q_0+1)}M_{\mu} } { V_{\mu} } \nonumber \\
&   \lesssim &  \mu^{-(\theta+1)},
\end{eqnarray}
where $\theta=r-q_0>0$. In light of  (\ref{Bb'}), we obtain the differential inequality 
\begin{eqnarray}\label{Bb3}
| \ln B(\mu) '|&  \lesssim &  \mu^{-\theta-1},
\end{eqnarray}
which by integration yields

\begin{eqnarray}\label{Bb4}
| \ln B(\mu) |&   \lesssim &   \mu^{-\theta} +C_0,
\end{eqnarray}
where $C_0>0$ is a constant independent of $\mu$.
Hence
\begin{eqnarray}\label{Aa5}
0< B(\mu) \lesssim C_1 e^{-\theta\mu^{-\theta}}  \sim  1.
\end{eqnarray}
which yields the uniform upper bound  for  $\mu\in [\mu_0, \infty)$,
The proof of Lemma \ref{lem8.1}  is now complete. 
\end{proof}

\begin{lemma}\label{lem8.2}
Suppose $(f_3)$ and $(f_4)$. Then, for all $\mu \in\left[\mu_0, \infty\right)$, we have 
$v_{\mu}(x)\lesssim|x|^{-N}$ for $|x| \geqslant R$ for some $R$ independent of $\mu$.
\end{lemma}

\begin{proof}
First, we claim that $v_{\mu} \in L^1$ holds for any fixed $\mu$. Indeed, recalling 
that $u_\lambda \in L^1$ holds by the uniform decay estimate in Lemma \ref{lem4.17},  
we conclude $v_{\mu} \in L^1$ by the definition of $v_\mu$ directly.

Integrating equation 
\begin{equation}
(-\Delta)^s v_{\mu}+v_{\mu} = \mu^{1-r}f(\mu v_{\mu})  
\end{equation}
over $\mathbb{R}^N$ and using the fact that $\int_{\mathbb{R}^N}(-\Delta)^s v_{\mu}=0$ holds, we obtain the identity
\begin{equation}\label{vf}
\int_{\mathbb{R}^N} v_{\mu}=\int_{\mathbb{R}^N} \mu^{1-r}f(\mu v_{\mu}) .
\end{equation}

Notice that
\begin{equation}
V_{\mu}=\mu^{-r} \int_{\mathbb{R}^N}  f(\mu v_{\mu})\mu v_{\mu}=
\int_{\mathbb{R}^N}v_{\mu}^r dx +  \mu^{-r}\int_{\mathbb{R}^N}g(\mu v_{\mu})\mu v_{\mu} dx,
\end{equation}
by $(f_4)$. Together with H\"older's inequality, this gives
\begin{eqnarray}\label{a}
\int_{\mathbb{R}^N} \mu^{1-r} f(\mu v_{\mu}) dx &\lesssim& \int_{\mathbb{R}^N} \mu^{1-r} (\mu v_{\mu})^{q-1}+ \mu^{1-r} (\mu v_{\mu})^{r-1}dx\nonumber \\
&\lesssim& \int_{\mathbb{R}^N} \mu^{q-r} ( v_{\mu})^{q-1}+ \int_{\mathbb{R}^N} ( v_{\mu})^{r-1}\nonumber \\
&\lesssim& \mu^{q-r} \left(\int_{\mathbb{R}^N} v_{\mu}\right)^{\frac{1}{q-1}}\left(\int_{\mathbb{R}^N} v_{\mu}^{q}\right)^{\frac{q-2}{q-1}}+   \left(\int_{\mathbb{R}^N} v_{\mu}\right)^{\frac{1}{r-1}}\left(\int_{\mathbb{R}^N} v_{\mu}^{r}\right)^{\frac{r-2}{r-1}}\nonumber \\
&\lesssim& \mu^{q-r} \left(\int_{\mathbb{R}^N} v_{\mu}\right)^{\frac{1}{q-1}}
+\left(\int_{\mathbb{R}^N} v_{\mu}\right)^{\frac{1}{r-1}},
\end{eqnarray}
where the last inequality follows from the uniform estimate proved in Lemma \ref{lem8.1}.
Going back to (\ref{vf}), we obtain
\begin{eqnarray}\label{a11}
\int_{\mathbb{R}^N} v_{\mu}
&\lesssim& \mu^{q-r} \left(\int_{\mathbb{R}^N} v_{\mu}\right)^{\frac{1}{q-1}}
+ \left(\int_{\mathbb{R}^N} v_{\mu}\right)^{\frac{1}{r-1}}.
\end{eqnarray}
Since $\mu\in [\mu_0,+\infty)$ and $q<r$, we deduce the uniform bound
\begin{eqnarray}\label{8.21}
\int_{\mathbb{R}^N}v_{\mu} \lesssim 1.
\end{eqnarray}
Since, the function $v_{\mu}(r)>0$ is decreasing in $r$. Hence, for any $|x|>0$,
\begin{equation}\label{8.22}
\int_{\mathbb{R}^N}v_{\mu} \geqslant \int_{|y| \leqslant |x|}v_{\mu}(|y|)dy \gtrsim 
|x|^N v_{\mu}(|x|) .
\end{equation}
In view of (\ref{8.21}) we conclude that $v_{\mu}(x)\lesssim|x|^{-N}$ holds for $|x| \geqslant R$ 
for some $R$ independent of $\mu$, completing the proof of Lemma \ref{lem8.2}.
\end{proof}

By $(f_4)$,
\[
\left(\frac{G(t)}{t^r}\right)' = O(t^{q-r-1}+t^{q_0 - r-1}),
\]
so there exists $\ell \ge 0$ such that
\[
\lim_{t \to \infty}\frac{G(t)}{t^r} = \ell, \quad \frac{f(t)}{t^{r-1}} \to a_\infty := 1+r\ell > 0.
\]

\begin{lemma} \label{lem8.3}
 Let $\{\mu_n\}_{n=1}^{\infty} \subset\left[\mu_0, \infty\right)$ be a sequence 
 such that $\mu_n \rightarrow \infty$ and suppose that $v_{\mu_n}(|x|)>0$ solves \eqref{v0} with $\mu = \mu_n$ and has Morse index $1$. 
 Then, after possibly passing to a subsequence, we have $v_{\mu_n} \rightarrow v_*$ in 
 $H^s(\mathbb{R}^N)$ as $n \rightarrow+\infty$. Moreover, the function $v_*(|x|)>0$ satisfies
\begin{equation}\label{8.23}
(-\Delta)^{s} v_*+v_*- a_\infty|v_*|^{r-2}v_* =0 \quad \text{in } \mathbb{R}^N,
\end{equation}
and has Morse index $1$.
\end{lemma}
\begin{proof}  Let $v_n=v_{\mu_n}$ in the following. Recall the uniform bound 
  $\left\|v_n\right\|_{H^{s}} \lesssim 1$ by the a priori bounds in Lemma \ref{lem8.1}. 
  Passing to a subsequence if necessary, we have $v_n \rightharpoonup v_*$ weakly in $H^{s}$. 
  Furthermore, by local Rellich compactness, we can assume that $v_n \rightarrow v_*$ in 
  $L_{\text {loc }}^p$, $2 \leqslant p < 2_s^*$ and pointwise a.e. in $\mathbb{R}^N$. 
  Thanks to the uniform decay bound in Lemma \ref{lem8.2}, we can upgrade this to 
  $v_n \rightarrow v_*$ strongly in $L^p$, $2 \leqslant p < 2_s^*$. Bootstrap arguments shows that 
  $v_n \rightarrow v_*$ in $H^{s}$. Obviously, $v_*\geqslant 0$ is nonzero and satisfies equation 
  (\ref{8.23}). Thus $0$ is the first eigenvalue of $L_{-, \infty} := 
  (-\Delta)^{s}+1-  a_\infty|v_*|^{r-2}$, showing $v_*>0$. Similar to the Morse index arguments in the proof of Lemma \ref{lem4.18}, we have $\mu(v_*)=1$. 
\end{proof}

\begin{proof}[Proof of Proposition \ref{propub}]
  Let $v_1, v_2 \in C^1([\mu_0,\infty),H_{\text{rad}}^{s})$ 
  be two maps such that for $i = 1,2$, $v_i(\mu) > 0$ solves \eqref{v0} and has Morse index $1$. For convenience, 
  we let $\kappa = 1/\mu$ and $\kappa_0 = 1/\mu_0$. Then for $i =1,2$, we consider maps 
  $v_i(\kappa) \in C^1((0,\kappa_0],H_{\text{rad}}^{s})$, and aim to show that 
  $v_1(\kappa) = v_2(\kappa)$ for all $\kappa \in (0,\kappa_0]$. By Lemma \ref{lem8.3}, 
  as $\kappa \to 0$, both $v_1(\kappa)$ and $v_2(\kappa)$ converge to $v_*$ which is positive, solves \eqref{8.23} and has Morse index $1$. By \cite{FS}, $v_*$ is unique (up to translations) and non-degenerate. 
  Then, by implicit function arguments like the proof of Lemma \ref{prop4.15}, we obtain that 
  $v_1(\kappa)=v_2(\kappa)$ for small positive $\kappa$, and thus for all $\kappa \in (0,\kappa_0]$. 
  The proof is complete.
\end{proof}

\subsection{Proof of Theorem \ref{thm} and Theorem \ref{thmN1}} \label{secpfmr}
\begin{proof}
It is readily seen that both cases satisfy the condition of Proposition \ref{propgb} and Proposition \ref{propub}.
Suppose on the contrary that \eqref{Q1} has two different positive ground state solutions $u_{1,\lambda_0}$ and 
$u_{2,\lambda_0}$ for some $\lambda = \lambda_0 > 0$. Then by Proposition \ref{propgb}, 
there exist two maps $u_1, u_2 \in C^1([\lambda_0,\infty),H_{\text{rad}}^{s})$ such that for $i = 1,2$, 
$u_i(\lambda)$ is positive, solves \eqref{Q1} and has Morse index $1$, with $u_i(\lambda_0) = u_{i,\lambda_0}$. Then, 
since $T_{\lambda}$ is one-to-one, Proposition \ref{propub} shows us that $u_1(\lambda) \equiv 
u_2(\lambda)$, a contradiction.
\end{proof}

\section{Uniqueness of the positive ground state solution: proof of Theorem \ref{thm:log}}
\label{sec:main-proofs}\label{sec:log-cont}
In this section we complete the proof of Theorem \ref{thm:log}. For the power mixtures we use the $\lambda$-continuation in Section~\ref{uniqueness}: the nonlinearity and $s$ remain fixed, and the large-frequency rescaling leads to the highest pure power. Here we use a nonlinearity deformation with both $0 < s < 1$ and $\lambda > 0$ fixed.

\subsection{The global branch corresponding to a nonlinearity deformation}
Write $P=p+\delta$ and for $0\le\theta\le1$, define
\begin{equation}\label{eq:log-path}
 f_\theta(t)=t^{P(1-\theta)}
       \bigl(ct^p\log(1+t^\delta)\bigr)^\theta
 =c^\theta t^{p+\delta(1-\theta)}[\log(1+t^\delta)]^\theta ~~ (t \ge 0); \quad f_\theta(t)= 0 ~~ (t < 0).
\end{equation}
Clearly, $f_0(t)=t^P$ and $f_1(t)=f(t)$. Note that $f_\theta$ satisfies $(f_1)$, $(f_2)$. Next, we check the conditions in Proposition \ref{thm:general-quotient-monotonicity}. A direct computation yields
\begin{align*}
\eta_\theta(t):=\frac{t f_\theta^{\prime}(t)}{f_\theta(t)}-1=p-1+\delta(1-\theta)+\frac{\theta \delta t^\delta}{\left(1+t^\delta\right) \log \left(1+t^\delta\right)} >0, \quad t > 0.
\end{align*}
In addition,
$$
t \eta_\theta^{\prime}(t)=\frac{\theta \delta^2 t^\delta\left[\log \left(1+t^\delta\right)-t^\delta\right]}{\left(1+t^\delta\right)^2\left[\log \left(1+t^\delta\right)\right]^2} \leq 0, \quad t > 0.
$$
Moreover, for $t > 0$ we have
$D^2\left(\frac{f_\theta}{t}\right)=\frac{f_\theta}{t}\left(\eta_\theta^2+t \eta_\theta^{\prime}\right),$
and
\begin{align*}
\eta_\theta^2+t \eta_\theta^{\prime}= & (p-1)^2 +2(p-1) \delta\left[1-\theta+\frac{\theta t^\delta}{\left(1+t^\delta\right) \log \left(1+t^\delta\right)}\right] \\
& +\delta^2(1-\theta)\left[1-\theta+\frac{2 \theta t^\delta}{\left(1+t^\delta\right) \log \left(1+t^\delta\right)}-\frac{\theta t^{2 \delta}}{\left(1+t^\delta\right)^2\left[\log \left(1+t^\delta\right)\right]^2}\right] \\
& +\frac{\theta \delta^2 t^\delta}{\left(1+t^\delta\right)^2 \log \left(1+t^\delta\right)} .
\end{align*}
Since $0<\frac{t^\delta}{\left(1+t^\delta\right) \log \left(1+t^\delta\right)}<1,$ we get
$\eta_\theta^2+t \eta_\theta^{\prime}>0$
for every $t>0$ and every $\theta \in[0,1]$.
Consequently,
$$
D^2\left(\frac{f_\theta(t)}{t}\right)>0, \quad t > 0.
$$
Thus $H_{f_\theta,\lambda}'(t) > 0$ for every $t > 0$.

\begin{lemma}\label{lem:log-peak}
Positive solutions of
$(-\Delta)^su+\lambda u=f_\theta(u)$, $0\le\theta\le1$ are
uniformly bounded in $L^\infty$.
\end{lemma}

\begin{proof}
Suppose by contradiction that there are parameters $\theta_n \in[0,1]$ and solutions $u_n$ such that $\theta_n \to \theta_*$ and
$M_n :=|u_n|_{\infty} \to \infty.$

Translate a maximum to the origin and set $\kappa_n=f_{\theta_n}\left(M_n\right)/M_n .$
Then  $\kappa_n \to\infty$. We define
\[
 w_n(x) :=M_n^{-1}u_n(\kappa_n^{-1/(2s)}x).
\]
Note that $0 < w_n(x) \leq 1 , w_n(0)=1$ and  $(-\Delta)^s w_n+\frac{\lambda}{\kappa_n} w_n=\frac{f_{\theta_n}\left(M_n w_n\right)}{f_{\theta_n}\left(M_n\right)} .$

For $\tau \in[0,1]$, define 

$$
\tilde{f}_n(\tau):=\frac{f_{\theta_n}\left(M_n \tau\right)}{f_{\theta_n}\left(M_n\right)} .
$$
For every fixed $\tau>0$, we have

$$
\tilde{f}_n(\tau)=\tau^{p+\delta\left(1-\theta_n\right)}\left[\frac{\log \left(1+M_n^\delta \tau^\delta\right)}{\log \left(1+M_n^\delta\right)}\right]^{\theta_n} \to \tau^{p+\delta\left(1-\theta_*\right)} .
$$
 Since each $\tilde{f}_n$ is nondecreasing and the limiting function is continuous, the convergence is uniform on $[0,1]$.
Consequently, if $w_n \rightarrow w$ locally uniformly, then $\tilde{f}_n\left(w_n\right) \to w^{p+\delta\left(1-\theta_*\right)} $ locally uniformly in $\mathbb{R}^N .$

 Interior compactness as in
Proposition~\ref{prop:peak} gives a bounded positive solution
$(-\Delta)^sw=w^{p+\delta\left(1-\theta_*\right)}$, $w(0)=1$.  A contradiction follows by the argument used in Proposition~\ref{prop:peak}.
\end{proof}

Consider:
\begin{align*}
  G: \ & H^s_{\text {rad }}(\mathbb{R}^N) \times [0,1]\to H^{-s}(\mathbb{R}^N), \\
  &G(u,\theta)=(-\Delta)^{s} u+ \lambda u -f_\theta(u).
\end{align*}
Note that $G(u,\theta)=0$ if and only if  $u \in H^s_{\text {rad }}$ solves
\begin{equation}\label{log*}
	(-\Delta)^s u+\lambda u=f_\theta(u) \quad \text { in } \mathbb{R}^N .
\end{equation}
In this subsection, 
we set
$$
L_{+,\theta} := (-\Delta)^s + \lambda - f'_\theta(u).
$$

\begin{lemma}\label{log4.15} 
 Fix $0 < \theta_0 \le 1$ and let $u_0 \in H_{\text {rad }}^{s}(\mathbb{R}^N)$ be a 
positive solution of $G(u_0, \theta_0)=0$ having Morse index $1$.   
Then for some $\delta>0$, there exists a map 
$u(\theta) \in C^1(  (\theta_0-\delta,\theta_0], H_{\text{rad}}^{s})$ 
such that the following assertions hold,  where we denote $u_\theta=u(\theta)$ in the sequel.

$(i)$ $u_\theta$ solves $G(u, \theta)=0$ and has Morse index $1$ for all $\theta \in (\theta_0-\delta,\theta_0]$.

$(ii)$ There exists $\epsilon>0$ such that $u_\theta$ is the unique solution of 
$G(u, \theta)=0$ for $\theta \in (\theta_0-\delta,\theta_0]$ in the 
neighborhood $\left\{u \in H_{\text {rad }}^{s}:\left\|u-u_\theta
\right\|_{H^{s}}<\epsilon\right\}$. In particular, we have that $u_{\theta_0}=u_0$ holds.
\end{lemma}

\begin{proof} 
We use the implicit function argument for the map
$$
\Psi: H_{\text {rad }}^{s} \times  (\theta_0-\delta,\theta_0]\rightarrow 
H_{\text {rad }}^{s}, \quad \Psi(u, \theta)=u-\frac{1}{(-\Delta)^s+\lambda} f_\theta( u).
$$
It is not difficult to see that $\Psi$ is a well-defined map of class $C^1$. 
Furthermore, $\Psi\left(u_0, \theta_0\right)=0$ and $\Psi(u,\theta)=0$ if 
and only if $u \in H_{\text {rad }}^{s}$ solves $G(u,\theta)=0$. Next, we consider the Fréchet derivative
$$
\partial _u \Psi\left(u_0, \theta_0\right)=1-K, K=\frac{1}{(-\Delta)^s+\lambda} f'_{\theta_0}(u_0)
$$
By  Theorem \ref{them2}, $0$ is not an eigenvalue of $L_{+,\theta_0}$ constrained on $L_{rad }^2$. Note that the 
operator $K$ is compact on $L_{\text {rad }}^2$ and we have that 
$1 \notin \sigma(K)$. Since $f'_{\theta_0}(u_0) \varphi \in H^{-s}$ if $\varphi \in 
H_{\text {rad }}^{s}$, we have that $K$ maps $H_{\text {rad }}^{s}$ to 
$H_{\text {rad }}^{s}$. Thus the bounded inverse $(1-K)^{-1}$ exists on 
$H_{\text {rad }}^{s}$, i.e., the Fréchet derivative $\partial_u\Psi$ has a 
bounded inverse on $H_{\text {rad }}^{s}$ at $\left(u_0, \theta_0 \right)$. 
By the implicit function theorem, we deduce that the assertions 
$(i)$ and $(i i)$ hold for some $\delta>0$ and 
$\epsilon>0$ sufficiently small.
\end{proof}

Next, we extend this branch. For this, we define
\begin{align*}
	\theta_*= \inf \{\hat \theta \ge 0: \ & u_{\theta}\in C^1((\hat\theta,\theta_0],
	H^s_{\text {rad }}(\mathbb{R}^N)), \\
	& u_{\theta}\text{ satisfies the assumptions of Lemma \ref{log4.15}}\text{ for } \theta\in (\hat \theta,\theta_0]\},
\end{align*}
and aim to prove $\theta_* = 0$. Similar to Lemma \ref{lem4.17} $(i)$, we know $u_{\theta} > 0$ for every $\theta\in ( \theta_*,\theta_0]$.

\begin{lemma}\label{theta81} 
	$u_\theta$ is uniformly bounded in $H^s$ with respect to $\theta \in (\theta_*,\theta_0]$.
\end{lemma}

\begin{proof}
	Lemma \ref{lem:log-peak} yields that
	$0<u_\theta(x) \leq M_*$ for some $M_* < \infty$. Recall that
\begin{align*}
\eta_\theta(t):=\frac{t f_\theta^{\prime}(t)}{f_\theta(t)}-1=p-1+\delta(1-\theta)+\frac{\theta \delta t^\delta}{\left(1+t^\delta\right) \log \left(1+t^\delta\right)} >0 .
\end{align*}
Thus there exists $\varepsilon>0$, independent of $\theta$, such that
$\varepsilon\leq \eta_\theta(t)\leq P-1$ for $0<t<M_*$ and $\theta\in[0,1]$. 
Hence,
\begin{equation}\label{log8.5}
(1+\varepsilon) f_\theta(t) \leq t f_\theta^{\prime}(t) \leq(p+\delta) f_\theta(t), \quad 0<t \leq M_* .
\end{equation}
Integrating (\ref{log8.5}) yields
\begin{equation}\label{log3}
(2+\varepsilon) F_\theta(t) \leq t f_\theta(t) \leq(p+\delta+1) F_\theta(t), \quad 0 < t \leq M_* .
\end{equation}
Viewing $f_\theta(t)$ as a function of  $\theta$ and differentiating gives
$$
\partial_\theta f_\theta(t)=f_\theta(t) \log \left(\frac{c \log \left(1+t^\delta\right)}{t^\delta}\right) .
$$
Since
$$
\lim _{t \to 0} \log \left(\frac{c \log \left(1+t^\delta\right)}{t^\delta}\right)=\log c,
$$
the logarithmic factor extends continuously to $t=0$ and is bounded on $\left[0, M_*\right]$. Thus there exists $K>0$, independent of $\theta$, such that

$$
\left|\partial_\theta f_\theta(t)\right| \leq K f_\theta(t), \quad 0 \leq t \leq M_* .
$$
Integrating in $t$, we obtain
\begin{align}\label{log1}
\left|\partial_\theta F_\theta(t)\right| & =\left|\int_0^t \partial_\theta f_\theta(a) d a\right| \leq K \int_0^t f_\theta(a) d a=K F_\theta(t).
\end{align}
Define
$J_\theta(u)=\frac12\|u\|_\lambda^2-\int_{\R^N}F_\theta(u)\,dx,$
where $F_\theta(t) = \int_0^tf_\theta(\tau)d\tau$ and $\|u\|_\lambda^2= \lambda|u|_2^2+\mathcal{E}(u,u).$
Clearly,
$$\frac{d}{d \theta} J_\theta\left(u_\theta\right)=\partial_\theta J_\theta\left(u_\theta\right)+D_u J_\theta\left(u_\theta\right)\left[\partial_\theta u_\theta\right] =\partial_\theta J_\theta\left(u_\theta\right).$$
Notice
$
\left\|u_\theta\right\|_\lambda^2=\int_{\mathbb{R}^N} f_\theta\left(u_\theta\right) u_\theta d x.$
Therefore, by (\ref{log3}) we have
\begin{equation}\label{log5}
\int_{\mathbb{R}^N} F_\theta\left(u_\theta\right) d x \leq \frac{1}{2+\varepsilon}\left\|u_\theta\right\|_\lambda^2 .
\end{equation}
Substituting  (\ref{log5}) into $J_\theta$, we find
\begin{align}\label{log6}
J_\theta\left(u_\theta\right) & \geq\left(\frac{1}{2}-\frac{1}{2+\varepsilon}\right)\left\|u_\theta\right\|_\lambda^2 =\frac{\varepsilon}{2(2+\varepsilon)}\left\|u_\theta\right\|_\lambda^2 .
\end{align}
In particular, $J_\theta\left(u_\theta\right)>0$, since $u_\theta \neq 0$. Combining (\ref{log1}), (\ref{log5}), and (\ref{log6}), we obtain 
\begin{equation}\label{eq:energy-gronwall}
 \left|\frac d{d\theta}J_\theta(u_\theta)\right|
 =\left|\int\partial_\theta F_\theta(u_\theta)\,dx\right|
 \le \frac{K}{2+\varepsilon}\|u_\theta\|_\lambda^2\le \frac{2K}{\varepsilon} J_\theta(u_\theta).
\end{equation}
 Hence, for any $ \theta \in\left(\theta_*, \theta_0\right], $
$$
J_\theta\left(u_\theta\right) \leq J_{\theta_0}\left(u_{\theta_0}\right) \exp \left(\frac{2 K}{\varepsilon}\left|\theta-\theta_0\right|\right) \leq J_{\theta_0}\left(u_{\theta_0}\right) e^{2 K / \varepsilon} .$$
Finally, (\ref{log6}) implies

$$
\left\|u_\theta\right\|_\lambda^2 \leq \frac{2(2+\varepsilon)}{\varepsilon} J_{\theta_0}\left(u_{\theta_0}\right) e^{2 K / \varepsilon}, \quad \theta \in\left(\theta_*, \theta_0\right] .
$$
This completes the proof.
\end{proof}

\begin{lemma}\label{decay}
There exists $C>0$ independent of $\theta \in (\theta_*, \theta_0]$ such that
$$
\|u_\theta\|_{L^1} \leq C \quad \text{and}\quad
0<u_\theta(x) \leq \frac{C}{1+|x|^N} .
$$
\end{lemma}

\begin{proof} 
We first use the uniform $L^2$ bound to show that the solutions tend to zero at infinity uniformly in $\theta$. By radial monotonicity of $u_\theta$, we have
$$
\left|B_r\right| u_\theta(r)^2 \leq \int_{B_r} u_\theta(x)^2 d x \leq \|u_\theta\|_{L^2}^2.
$$
Consequently,
\begin{equation}\label{re1}
u_\theta(r) \leq C r^{-N / 2}, \quad r>0,
\end{equation}
with $C$ independent of $\theta$.
As $t \to 0$, we have
$$
\frac{f_\theta(t)}{t} \leq \max \{1, c\} t^{p+\delta-1} \to 0
$$
uniformly in $\theta$. We may therefore choose $a>0$ such that

$$
f_\theta(t) \leq \frac{\lambda}{2} t, \quad 0 \leq t \leq a,
$$
for every $\theta \in[0,1]$. By (\ref{re1}), there exists $R_0>0$, independent of $\theta$, such that $u_\theta(x) \leq a$ whenever $|x| \geq R_0$. Hence
\begin{equation}\label{re2}
f_\theta\left(u_\theta(x)\right) \leq \frac{\lambda}{2} u_\theta(x), \quad|x| \geq R_0 .
\end{equation}
Rewrite the equation (\ref{log*}) as
\begin{equation}\label{re3}
\left((-\Delta)^s+\frac{\lambda}{2}\right) u_\theta=f_\theta\left(u_\theta\right)-\frac{\lambda}{2} u_\theta .
\end{equation}
We use the following expression
$$
u_\theta=G_{s, \lambda / 2} *\left(f_\theta\left(u_\theta\right)-\frac{\lambda}{2} u_\theta\right),
$$
where $ G_{s, \lambda / 2} $ denotes the kernel of 
$\left((-\Delta)^s+\frac{\lambda}{2}\right)^{-1}$. This kernel is nonnegative and belongs to $L^1\left(\mathbb{R}^N\right)$, with 
\begin{equation}\label{re4}
\int_{\mathbb{R}^N} G_{s, \lambda / 2}(x) d x=\frac{2}{\lambda} .
\end{equation}
By positivity of the kernel, we have
\begin{equation}\label{re5}
 u_\theta \leq G_{s, \lambda / 2} *\left(f_\theta\left(u_\theta\right)-\frac{\lambda}{2} u_\theta\right)^{+},
\end{equation}
where $v^{+}=\max \{v, 0\}$.
By (\ref{re2}), the positive part on the right-hand side vanishes outside $B_{R_0}$. Inside this ball, it satisfies

$$
\left(f_\theta\left(u_\theta\right)-\frac{\lambda}{2} u_\theta\right)^{+} \leq f_\theta\left(u_\theta\right) \leq \max \{1, c\} M_*^{p+\delta} .
$$
Therefore,
\begin{equation}\label{re6}
\int_{\mathbb{R}^N}\left(f_\theta\left(u_\theta\right)-\frac{\lambda}{2} u_\theta\right)^{+} d x \leq \max \{1, c\} M_*^{p+\delta}\left|B_{R_0}\right| .
\end{equation}
Integrating (\ref{re5}) and  applying (\ref{re4}) and (\ref{re6}), we obtain

$$
\begin{aligned}
\int_{\mathbb{R}^N} u_\theta(x) d x & \leq\left(\int_{\mathbb{R}^N} G_{s, \lambda / 2}(x) d x\right) \int_{\mathbb{R}^N}\left(f_\theta\left(u_\theta\right)-\frac{\lambda}{2} u_\theta\right)^{+} d x \\
& \leq \frac{2}{\lambda} \max \{1, c\} M_*^{p+\delta}\left|B_{R_0}\right|.
\end{aligned}
$$
Using
$
\left|B_r\right| u_\theta(r) \leq \int_{B_r} u_\theta(x) d x$ we get
$u_\theta(r) \leq C r^{-N}$ for some $C  > 0$ independent of $\theta$.
Combining this estimate with $u_\theta \leq M_*$, we conclude that
$$
u_\theta(x) \leq \frac{C}{1+|x|^N}.
$$
 The proof is complete.
\end{proof}

\begin{lemma}\label{positive}
 Let $\{\theta_n\}_{n=1}^{+\infty}\subset 
(\theta_*,\theta_0]$ be a sequence such that $\theta_n\to \theta_*$. Then, after possibly passing to a 
subsequence, we have $u_{\theta_n}\to u_{\theta_*}$ in 
$H^s_{\text {rad }}(\mathbb{R}^N)$, where $u_{\theta_*}(|x|)>0$ is a positive solution of the following equation
\begin{equation}\label{e4.34}
(-\Delta)^s u_{\theta_*}+\lambda u_{\theta_*}=f_{\theta_*}\left(u_{\theta_*}\right) \quad \text{in } \mathbb{R}^N,
\end{equation}
having Morse index $1$.
\end{lemma}

\begin{proof}
By Lemma \ref{theta81}, $\{u_{\theta_n}\}$ is bounded in $H^s$. Passing to a 
subsequence if necessary, we assume that $u_{\theta_n} \rightharpoonup  u_{\theta_*}$ 
in $H^s, u_{\theta_n} \rightarrow u_{\theta_*}$ in $L_{loc}^r$ for 
$1 < r<2_s^*, u_{\theta_n} \rightarrow u_{\theta_*}$ a.e. on $\mathbb{R}^N$. 
With the help of the uniform decay estimate shown in Lemma \ref{decay}, 
we obtain that $u_{\theta_n} \rightarrow u_{\theta_*}$ in $L^r$ for 
$1< r<2_s^*$. Bootstrap arguments implies that $u_{\theta_n} \rightarrow 
u_{\theta_*}$ in $H^{s}$. One can check that $u_{\theta_*}\geqslant 0$ is nonzero and satisfies equation 
(\ref{e4.34}). Hence, $u_{\theta_*}=\left((-\Delta)^{s}+\lambda\right)^{-1} f_{\theta_*}\left(u_{\theta_*}\right)$. Since the kernel of $\left((-\Delta)^{s}+\lambda\right)^{-1}$ is positive, we obtain that $u_{\theta_*}>0$ as well. Similar to the Morse index arguments in proof of Lemma \ref{lem4.18}, we have $\mu(u_{\theta_*})=1$.

\end{proof}

\begin{lemma}\label{branch}
 We conclude that $\theta_*=0$.
\end{lemma}

\begin{proof}
  Suppose on the contrary that $\theta_* >0$. By Lemma \ref{positive}, 
  we get $u_{\theta_*}$ which is a positive radial solution of \eqref{e4.34} with Morse index $1$. Applying Theorem \ref{them2}, 
  we obtain that $u_{\theta_*}$ is non-degenerate. Then using implicit function arguments again, we can extend the branch $u(\theta)$ to $\theta_* - \delta$, 
  contradicting the definition of $\theta_*$. Hence, we complete the proof.
\end{proof}

\subsection{Proof of Theorem \ref{thm:log}}
\begin{proof}[Proof of Theorem~\ref{thm:log}]
In view of Lemmas \ref{positive}, \ref{branch}, we can complete the proof by repeating discussions in Subsection \ref{secpfmr}.
\end{proof}

\medskip

\noindent \textbf{Acknowledgements:} Song is supported by the China Postdoctoral Science Foundation (2024T170452). This work was completed when Song was visiting the Universit\'e Marie et Louis Pasteur, LmB (UMR 6623).

\end{document}